\documentclass{amsart}
  
  \usepackage{amsthm} 
  \usepackage{amsmath} 
  \usepackage{amssymb} 
  \usepackage{mathtools} 
   \usepackage{amsfonts}
  \usepackage[dvipsnames]{xcolor}
  \usepackage[T1]{fontenc}
  
  \usepackage{wrapfig}
   
  \usepackage[scr=boondox]{mathalfa}

 \usepackage{graphicx}
  \usepackage{xspace}
  \usepackage[all]{xy}
  \usepackage{pinlabel}
  \usepackage{enumerate}
\usepackage[margin=1.3in]{geometry}
  \usepackage{hyperref}

  \newcommand{\calA}{\mathcal{A}}
  \newcommand{\calB}{\mathcal{B}}

  \newcommand{\calF}{\mathcal{F}}

  \newcommand{\calM}{\mathcal{M}}

  \newcommand{\calS}{\mathcal{S}}

  \newcommand{\NN}{\mathbb{N}}

  \newcommand{\RR}{\mathbb{R}}

  \newcommand{\ZZ}{\mathbb{Z}}

  \newtheorem{theorem}{Theorem}[section]
  \newtheorem{proposition}[theorem]{Proposition}
  
  \newtheorem{lemma}[theorem]{Lemma}

  \newtheorem{question}[theorem]{Question}

  \theoremstyle{definition}
  \newtheorem{definition}[theorem]{Definition}

  \newtheorem*{claim*}{Claim}

  \newtheorem{fact}[theorem]{Fact}
  \newtheorem*{question*}{Question}
  \newtheorem*{answer*}{Answer}
  \newtheorem*{application*}{Application}
  \newtheorem*{notation*}{Notation}

  \theoremstyle{remark}
  \newtheorem{remark}[theorem]{Remark}
  \newtheorem*{remark*}{Remark}
  
  \newcommand{\secref}[1]{Section~\ref{#1}}
  \newcommand{\thmref}[1]{Theorem~\ref{#1}}
  
  \newcommand{\lemref}[1]{Lemma~\ref{#1}}
  \newcommand{\propref}[1]{Proposition~\ref{#1}}

  \newcommand{\figref}[1]{Figure~\ref{#1}}
  \newcommand{\defref}[1]{Definition~\ref{#1}}
  
  \newcommand{\eqnref}[1]{Equation~\eqref{#1}}

  \DeclareMathOperator{\Mod}{Mod}

  \DeclareMathOperator{\codim}{codim}  

  \newcommand{\Map}{\ensuremath{\operatorname{Map}}\xspace}    
       
        \newcommand{\FMap}{\ensuremath{\operatorname{FMap}}\xspace}  
   
    \newcommand{\Accu}{\ensuremath{\mathrm{Accu}}\xspace}  

\newcommand{\Qinf}{{Q^{\infty}}}

  \newcommand{\param}{{\mathchoice{\mkern1mu\mbox{\raise2.2pt\hbox{$
  \centerdot$}}
  \mkern1mu}{\mkern1mu\mbox{\raise2.2pt\hbox{$\centerdot$}}\mkern1mu}{
  \mkern1.5mu\centerdot\mkern1.5mu}{\mkern1.5mu\centerdot\mkern1.5mu}}}

  \renewcommand{\setminus}{{\smallsetminus}}
  
  \newcommand{\st}{\mathbin{\mid}} 
   
  \newcommand{\from}{\colon\thinspace}

  \DeclarePairedDelimiter\abs{\lvert}{\rvert}
  \DeclarePairedDelimiter\norm{\lVert}{\rVert}

  \newcommand{\xdownarrow}[1]{%
  {\left\downarrow\vbox to #1{}\right.\kern-\nulldelimiterspace}}

\definecolor{cyan1}{RGB}{0,255,255}
\definecolor{dustyblue}{RGB}{55,171,200}

\begin{document}

\title{Asymptotic Dimension of Big Mapping Class Groups}


 \author   {Curtis Grant}
 \address{Department of Mathematics, University of Toronto, Toronto, ON }
 \email{curtis.grant@mail.utoronto.ca}

 \author   {Kasra Rafi}
 \address{Department of Mathematics, University of Toronto, Toronto, ON }
 \email{rafi@math.toronto.edu}

 \author   {Yvon Verberne}
 \address{Department of Mathematics, University of Toronto, Toronto, ON }
 \email{yvon.verberne@mail.utoronto.ca}
 
 
  \date{\today}

\begin{abstract}
Even though big mapping class groups are not countably generated, certain big mapping class groups can be 
generated by a coarsely bounded set and have a well defined quasi-isometry type. We show that the
big mapping class group of a stable surface of infinite type with a coarsely bounded generating set 
that contains an essential shift has infinite asymptotic dimension. This is in contrast with the mapping class 
groups of surfaces of finite type where the asymptotic dimension is always finite. We also give a topological 
characterization of essential shifts. 
\end{abstract}

\maketitle

\section{Introduction}

In this paper, the surface $\Sigma$ is an orientable, connected, second-countable 
$2$--manifold without boundary. We further assume that $\Sigma$ is stable, that is, every end of $\Sigma$ 
has a stable neighborhood (see \defref{Def:stable}). The mapping class group of $\Sigma$, denoted by 
$\Map(\Sigma)$, is the group of orientation preserving homeomorphisms of $\Sigma$ 
up to isotopy.  A surface $\Sigma$ is said to be of finite type when $\pi_{1}(\Sigma)$ 
is finitely generated and is of infinite type otherwise. The mapping class groups
of surfaces of infinite type are referred to as big mapping class groups. 

When $\Sigma$ is a surface of finite type, $\Map(\Sigma)$ is finitely generated. 
A finite generating set defines a word metric on $\Map(\Sigma)$ which is well-defined
up to quasi-isometry independent of the particular finite generating set. The large 
scale geometry of the mapping class group, that is the geometry of the quasi-isometry class of 
such metrics, has been studied extensively \cite{BKMM, Bowditch2, Bowditch1, EMR, Hamenstadt, MM2}.

In contrast, big mapping class groups are not even countably generated. However, using the 
framework of Rosendal for coarse geometry of non locally compact groups \cite{Rose},
we can establish a notion large scale geometry for big mapping class groups
when $\Map(\Sigma)$ has a coarsely bounded generating set. For a Polish topological group 
$G$, a subset $A \subset  G$ is coarsely bounded, abbreviated CB, if every compatible 
left-invariant metric on $G$ gives $A$ finite diameter. We say $G$ is locally CB
if some neighborhood of the identity in $G$ is CB, and we say $G$ is CB generated if 
$G$ has a generating set that is a union of a CB neighborhood of the identity and a finite set. 
Such a generating set defines a word metric on $G$ that is well defined up to a quasi-isometry.
Namely, word metrics on $G$ associated to different CB generating sets are quasi-isometric
to each other. 

Mann-Rafi gave a classification of mapping classes of stable surfaces that are CB generated \cite{MR}. 
We are interested in the study of the coarse geometry of such big mapping class groups. In particular, 
we would like to know if big mapping class groups have finite asymptotic dimension.

The notion of asymptotic dimension was introduced by Gromov in \cite{Gromov} where he also 
proved that $\delta$-hyperbolic groups have finite asymptotic dimension. Many other groups have also
been shown to have finite asymptotic dimensions \cite{BF, Ji, Osin, Roe}. The study of the asymptotic 
dimension in mapping class groups started with the work of Bell-Fujiwara \cite{BF} who modified Gromov's 
argument to show that the curve complex has finite asymptotic dimension. Masur and Minsky showed that the 
curve complex is Gromov hyperbolic \cite{MM1} and that the geometry of various curve complexes are closely 
linked with the coarse geometry of the mapping class group \cite{MM2}. Using these facts, Bestvina, 
Bromberg and Fujiwara showed that, for surfaces of finite type, $\Map(\Sigma)$ has finite asymptotic 
dimension by embedding $\Map(\Sigma)$ in a finite product of \emph{trees of curve complexes} \cite{BBF}. 

In comparison with surfaces of finite type, there are new phenomena present in 
big mapping class groups. For example, there are homeomorphisms where parts of the surface 
are pushed to infinity (never recurring back). The simplest form of such a map is a shift map. 
Consider an infinite strip in $\RR^2$ with $\ZZ$ acting on the strip by translations. Cut out an
equivariant family of disks and attach identical surfaces $\Sigma_i$ (possibly of infinite type) to the
boundaries of these disks. Then $\ZZ$ still acts by translations on the strip, sending the surface 
$\Sigma_i$ homeomorphically to the surface $\Sigma_{i+1}$. Embedding this in a larger strip $\sigma$, 
we can construct a homeomorphism $h_\sigma$ of $\sigma$ that acts as described above in the smaller 
strip, but the restriction of $h_\sigma$ to the boundary of $\sigma$ is the identity. Assume $\Sigma$
contains a copy of $\sigma$ where the ends of $\sigma$ exit different ends of $\Sigma$. Then there is a
homeomorphism of $\Sigma$ (again called $h_\sigma$) that is as described in $\sigma$ and is the identity 
map outside of $\sigma$.  We call $h_\sigma$ a \emph{shift map}. We say a shift map 
$h_\sigma \in \Map(\Sigma)$ is \emph{essential} if the group generated by $h_\sigma$, $\langle h_\sigma \rangle$, 
is not a coarsely bounded subgroup of $\Map(\Sigma)$.

\begin{theorem}[Main Theorem]\label{Thm:Main}
Assume $\Sigma$ is stable and $\Map(\Sigma)$ is CB generated. If $\Map(\Sigma)$ 
contains an essential shift, then the asymptotic dimension of $\Map(\Sigma)$ is infinite.
\end{theorem}

That is, the existence of an essential shift results in $\Map(\Sigma)$ having 	
very non-trivial geometry. In fact, any subset of $\ZZ^n$ for any $n>0$ can be embedded 
quasi-isometrically in $\Map(\Sigma)$ (see \thmref{Thm:Infinite-Asymptotic-Dimension-Qinf}). 

We also provide the topological classification of essential shifts. To do this, it is more natural to work with a 
certain finite index subgroup of $\Map(\Sigma)$.  Let $E$ be the end space of $\Sigma$. Mann-Rafi \cite{MR} 
defined a partial order on $E$ measuring the local complexity. Let $\FMap(\Sigma)$ be the subgroup of 
$\Map(\Sigma)$ that fixes the set of \emph{isolated maximal points} in the end space $E$ 
(see Section \ref{SubSec:FiniteIndexSubgroup} for the definition).  In the setting of CB (but infinite) generated 
groups, it is not always true that a finite index subgroup $H$ of $G$ is quasi-isometric to $G$. However, we show: 

\begin{theorem}
The group $\FMap(\Sigma)$ is quasi-isometric to $\Map(\Sigma)$. 
\end{theorem}

It turns out essential shifts are present when the end space is \emph{two-sided} in a certain sense.  
Let $E(z)$ denote the $\FMap(\Sigma)$ orbit of $z$ in $E$ and $\Accu(z)\subset E$ be the accumulation set 
of $E(z)$. Also let $E^G \subset E$ be the set of non-planar ends of $\Sigma$
(for $x \in E^G$, every neighborhood of $x$ in $\Sigma$ has non-zero genus). 

\begin{definition}[Two-sided] \label{Def:two-sided}
We say $E(z)$ is \emph{two-sided} if $E(z)$ is countable and $\Accu(z) = X \sqcup Y$ where 
$X, Y$ are non-empty disjoint closed $\FMap(\Sigma)$--invariant subsets of $E$. 
We say $E^G$ is \emph{two-sided} if $E^G = X \sqcup Y$ where 
$X, Y$ are non-empty disjoint closed  $\FMap(\Sigma)$--invariant subsets of $E$. 
We say $E$ is two-sided if either $E(z)$ is two-sided for some $z \in E$ or if $E^G$ is two-sided. 
\end{definition}

\begin{theorem}[Existence of essential shift]\label{Thm:Existence}
$\Map(\Sigma)$ contains an essential shift if and only if the end space $E$ of $\Sigma$ is two-sided. 
\end{theorem}

We can refine this theorem to give a characterization of exactly which shift maps are essential. For 
a shift map $h_\sigma$ with support $\sigma$, let $x, y \in E$ be the ends of $\Sigma$
associated to the ends of $\sigma$, that is, the ends $\sigma$ exits towards. Let $E(\Sigma_i )$
be set of ends of $\Sigma_i$, where the $\Sigma_i$ are the subsurfaces shifted by $h_{\sigma},$ and 
let $\calM(\Sigma_i)$ be the set of maximal points in $E(\Sigma_i)$ (see \secref{SubSec:partial-order}).

\begin{theorem}[Topological characterization of an essential shift]\label{Thm:Characterization}
A shift map $h_\sigma$ is essential if and only if either
\begin{itemize}
\item $\Sigma_i$ has finite genus and $E^G$ is two-sided giving a decomposition $E^G=X \sqcup Y$
where $x \in X$ and $y \in Y$; or 
\item for some $z \in \calM(\Sigma_i)$, $E(z)$ is two-sided giving a decomposition $\Accu(z)=X \sqcup Y$
where $x \in X$ and $y \in Y$.
\end{itemize}
\end{theorem}

There are other topological properties of $\Sigma$ that result in $\Map(\Sigma)$ having a non-trivial geometry.
A subsurface $R$ of $\Sigma$ is called non-displaceable if for every $f \in \Map(\Sigma)$, $f(R)$ intersects $R$. 
In \cite{HQR}, it was shown that $\Map(\Sigma)$ acts non-trivially on an infinite diameter hyperbolic space if and 
only if $\Sigma$ contains a non-displaceable subsurface. Using the classification of coarsely bounded mapping 
class groups given in \cite{MR}, we show:

\begin{theorem} [Two sources of non-trivial geometry] \label{Thm:Non-Trivial-Geometry}
If $\Map(\Sigma)$ does not have an essential shift and $\Sigma$ does not contain 
a non-displaceable subsurface then $\Map(\Sigma)$ is quasi-isometric to a point.
\end{theorem}

That is, any surface that does have an essential shift map or a non-displaceable subsurface does not
have interesting geometry. The finiteness of the asymptotic dimension of the mapping class groups in the remaining cases 
is still open. 

\begin{question}
Let $\Sigma$ be a tame infinite type surface that contains a non-displaceable subsurface such that $\Map(\Sigma)$ 
has no essential shifts. Is $\mathrm{asdim}(\Map(\Sigma))$ always finite? 
\end{question}

\begin{remark}
The classification of CB generated big mapping class groups in \cite{MR} was
carried out in a larger class of \emph{tame} surfaces where only certain 
ends are assumed to have stable neighborhoods. Most of the arguments in this paper
still work in the setting of tame surfaces as well, for example, an essential 
shift does always imply infinite asymptotic dimension. However, without the
assumption that $\Sigma$ is stable, $E$ being two sided does not imply 
existence of a shift map. Hence the statement and some arguments are cleaner 
with this assumption. The class of stable surfaces, first used in 
\cite{FGM}, is a large and natural class of surfaces to work with and it 
includes all easily constructed infinite type surfaces. 
\end{remark}

\subsection*{Outline of the paper}
We show $\FMap(\Sigma)$ has infinite asymptotic dimension by embedding 
the infinite dimensional cube, $Q^{\infty} \subset \{0,1\}^\NN$, quasi-isometrically into 
$\FMap(\Sigma)$. We introduce the definition of asymptotic dimension in 
\secref{Sec:AsymptoticDimension} and show that $\Qinf$ has infinite asymptotic dimension. 
In \secref{Sec:FiniteIndexSubgroup}, we introduce the 
finite-index subgroup of the mapping class group $\FMap(\Sigma)$ and show that 
it is quasi-isometric to the mapping class group, $\Map(\Sigma)$.

In Section \ref{Sec:SharkTank}, we carry out our main arguments in the simple
case of the \emph{shark tank}, $T$, which is a cylinder with a $\ZZ$ action 
and a discrete set of punctures exiting both ends of the cylinder. We define a length function on 
$\FMap(T)$ which counts how many punctures from one side of the cylinder 
are mapped to the other side. We then use this length function to show that 
there is a quasi-isometric embedding from $\Qinf$ to $\FMap(T)$. 

The proof in the general case is essentially the same, 
except one has to define appropriate length functions whenever $E$ is two-sided.
When $E(z)$ is two-sided for some $z \in E$, the length function counts 
the number of ends in $E(z)$ that are moved from one side to another. 
When $E^G$ is two-sided, we rely on the action of $\FMap(\Sigma)$ on homology 
to construct a suitable length function (see Section~\ref{Sec:Length-Function}). 
Using these length functions, we again show that there is a quasi-isometric
embedding of $\Qinf$ to $\FMap(\Sigma)$ whenever $E$ is two-sided which
implies $\FMap(\Sigma)$ has infinite asymptotic dimension (see \secref{Sec:InfAsympDimGeneral}). 

Theorems \ref{Thm:Existence} and \ref{Thm:Characterization} are proven 
in \secref{Sec:EquivalenceAlgebraicTopologicallyEssential}
and \thmref{Thm:Non-Trivial-Geometry} is proven in \secref{Sec:NonTrivialGeometry}.

\subsection*{Acknowledgements} We would like to thank Jonah Gaster for helpful conversations.
The first author was supported by an NSERC-USRA award.
The second author was supported by by an NSERC Discovery grant, RGPIN 06486.
The third author was supported by the National Science Foundation under Grant No. DMS-1928930 while participating 
in a program hosted by the Mathematical Sciences Research Institute in Berkeley, California, during the Fall 2020 semester. 
The second author was also partially supported by an NSERC-PDF Fellowship.

\medskip 
\section{Asymptotic Dimension}\label{Sec:AsymptoticDimension}

\subsection{Asymptotic dimension}\label{SubSec:AsymptoticDimension}
In this section, we review some facts regarding the asymptotic dimension. 
The notion of asymptotic dimension is due to Gromov 
\cite{Gromov}. For a full discussion of asymptotic dimension see \cite{BD}.

\begin{definition}[Asymptotic Dimension]\label{Def:AsymptoticDimension}
Let $X$ be a metric space. We say that $\mathrm{asdim}(X) \leq n$ if for every 
$R >0$ there exists a covering of $X$ by open sets $\{U_i\}_{i=1}^{\infty}$ such that 
$\sup_{i \in \mathbb{N}} \{\text{diam}(U_i)\} < \infty$,
and every ball of radius $R$ in $X$ intersects at most $n+1$ elements of the cover $\{U_i\}$. 
The asymptotic dimension is the least $n$ for which $\text{asdim}(X) \leq n$. We then write 
$\mathrm{asdim}(X) = n$. If no such $n$ exists, then we say $X$ has infinite asymptotic dimension.
\end{definition}

We now recall two useful facts. The first fact states high dimensional 
space cannot be coarsely embedded in a small dimensional space. 

\begin{fact}[Theorem 5 of \cite{BD}]\label{Fact:EmbeddingAsympDim}
If $f\from X \to Y$ is a quasi-isometric embedding between two metric spaces $X$ and $Y$, 
then $\text{asdim(X)} \leq \text{asdim}(Y)$.
\end{fact}

The second fact is regarding the asymptotic dimension of $\ZZ^{n}$. 

\begin{fact}[Theorems 5 and 6 in \cite{BD}]\label{Fact:AsympDimZ}
The asymptotic dimension of $\ZZ^n$ with the metric 
\[
d_{\mathbb{Z}^n}\big( (a_1,..,a_n),(b_1,..,b_n)\big) = \sum_{i=1}^{n} \abs{a_i-b_i},
\]
is exactly $n$.
\end{fact}

\subsection{The space $\Qinf$}\label{SubSec:Qinf}
As mentioned in the introduction, we prove $\FMap(\Sigma)$ has infinite asymptotic 
dimension by embedding a copy of an infinite cube $\Qinf$ into $\FMap(\Sigma)$. 

\begin{definition}[$\Qinf$]\label{Def:Qinf}
Let $\Qinf \subset \{0,1\}^{\NN}$ be the set of sequences that are eventually $0$. 
We equip $\Qinf$ with the $\ell^1$ metric, so that two sequences $a=(a_1,a_2,a_3,...)$ 
and $b=(b_1,b_2,b_3,...)$ have distance 
\[
d_{\Qinf}(a,b) = \sum_{i=1}^{\infty} \abs{a_i-b_i}. 
\]
This sum converges since all but finitely many terms are equal to $0$. 
We denote the zero sequence simply by $0$ and denote 
$d_{\Qinf}(0,a)$ by $|a|$. Thinking of $\{0, 1\}$ as the group $\ZZ/2\ZZ$, 
we can then write $\abs{a-b}$ to denote $d_{\Qinf}(a,b)$. 
\end{definition}

\begin{theorem}\label{Thm:Infinite-Asymptotic-Dimension-Qinf} 
The space $\Qinf$ has infinite asymptotic dimension. In fact, for every integer $n\geq1$, 
$\ZZ^n$ quasi-isometrically embeds in $\Qinf$ 
\end{theorem}

\begin{proof}
We begin by proving the second assertion. For every odd prime number $p \in \NN$, consider 
the map $f_p \from \ZZ \to \Qinf$ defined pointwise as follows: for a positive integer $m>0$ define
\begin{align*}
    f_p(m)_{j} = 
    \begin{cases}
    1 &\text{if} \ j=p^{k} \ \text{and} \ 0< k \leq m
    \\
    0 &\text{otherwise}
    \end{cases} 
\end{align*}
for a negative integer $m<0$ define
\begin{align*}
    f_p(m)_{j} = 
    \begin{cases}
    1 &\text{if} \ j=2p^{k} \ \text{and} \ 0< k \leq \abs{m}
    \\
    0 &\text{otherwise}
    \end{cases}
\end{align*}
and for $m=0$ we define $f_p(m)$ to be the zero sequence. We observe that the map $f_p$ is an isometry 
from $\ZZ$ to $\Qinf$. In fact, for $m_1, m_2 \in \ZZ$, if $m_1$ and $m_2$ have the same sign, 
then $d_{\Qinf}(f_p(m_1),f_p(m_2))$ is the number of positions in which $f_p(m_1)$ and $f_p(m_2)$ differ, 
which is precisely $\abs{m_1-m_2} = d_{\ZZ}(m_1,m_2)$ and if $m_2 \leq 0 \leq m_1$ then 
$d_{\Qinf}(f_p(m_1),f_p(m_2))$ is the number of positions where $f_p(m_1)$ is $1$ plus 
the number of positions where $f_p(m_2)$ is one which is $m_1 + \abs{m_2} = d_{\ZZ}(m_1,m_2)$. 

For $n\geq 1$, we choose $n$ distinct odd primes $p_1,...,p_n$. Now we define 
$F \from \ZZ^n  \to \Qinf$ as follows. Given a tuple ${\mathbf{ a}}=(a_1,...,a_n) \in \ZZ^n$ let 
\[
 F({\mathbf{a}}) = \sum_{i=1}^{n} f_{p_i}(a_i).
\]
This map is one-to-one since, for different values of $i$, the indices on which each $f_{p_i}$ are 
non-zero are distinct. Furthermore, 
\[
 \abs{F({\mathbf{a}})} =  \sum_{i=1}^{n} a_i.
\]
Therefore,  $F$ is an isometry which proves the second assertion.  

As stated in Fact \ref{Fact:AsympDimZ}, the asymptotic dimension of $\ZZ^n$ equipped with the $\ell^1$ 
metric is $n$.  Hence, the second assertion of the Theorem and Fact \ref{Fact:EmbeddingAsympDim} imply that 
$\text{asdim}(\Qinf) \geq n$ for every $n$, which proves the first assertion.
\end{proof}

\section{A finite index subgroup of the mapping class group}\label{Sec:FiniteIndexSubgroup}

Recall from the introduction that $\Sigma$ is an infinite-type, orientable, connected, 
second-countable $2$--manifold without boundary with tame end space such that 
$\Map(\Sigma)$ has a CB generating set. In this section, we recall the definition of tame
end space and describe the CB generating set for $\Map(\Sigma)$. We also introduce the 
finite-index subgroup $\FMap(\Sigma)$ of $\Map(\Sigma)$ and show that 
$\FMap(\Sigma)$ is quasi-isometric to $\Map(\Sigma)$. Then
Fact \ref{SubSec:AsymptoticDimension} implies that these groups have the same 
asymptotic dimension. 

\subsection{A partial order on the space of ends of $\Sigma$}\label{SubSec:partial-order}
Each topological $X$ has an \textit{space of ends} which is defined as 
the inverse limit $\varprojlim_{K\subset \Sigma} \pi_0(X \setminus K)$ 
as K ranges over the compact subsets of $X$. We denote the space of
end of $\Sigma$ by $E=E(\Sigma)$. For any subsurface $\Sigma' \subset \Sigma$,
$E(\Sigma')$ is the space of end of $\Sigma'$. We always
assume subsurfaces of $\Sigma$ have compact boundary hence, to ensure
$E(\Sigma') \subset E(\Sigma)$. 

A point $x \in E$ is a \emph{non-planar end} or is \emph{accumulated by genus}
if every neighborhood of $x$ in $\Sigma$ has non-zero genus. Otherwise, 
$x$ is a \emph{planar end} and it admits a neighborhood which can be 
embedded in the plane. We denote the subset of non-planar ends of $\Sigma$ 
by $E^G$. Topologically, $E$ is closed and 
totally disconnected and hence it is homeomorphic to a closed subset of 
the Cantor set. The space $E^G$ is a closed subset of $E$.

Richards proved that orientable, boundary-less, infinite type surfaces are completely 
classified by their genus (possibly infinite), the space of ends $E$, and the subset of 
ends accumulated by genus $E^{G}$ \cite{Rich}.  When we talk about homeomorphisms between 
subsets of $E$, we always assume they are \emph{type preserving}. That is, we say $U \subset E$ 
is homeomorphic to $V \subset E$ if there is a homeomorphism $f \from U \to V$ such that $f$ 
sends $U \cap E^G$ homeomorphically to $V \cap E^G$. Every homeomorphism of $\Sigma$ induces a 
(type preserving) homeomorphism on the space of ends and, by Richards classification, every 
homeomorphism of $E$ is induced by some element of $\Map(\Sigma)$.  

The following definition, given by Mann and Rafi, gives a ranking of the \emph{local complexity of an end} providing 
a partial order on equivalence classes of ends. See \cite[Section 4 and Section 6.3]{MR} for a detailed discussion.

\begin{definition}\label{Def:BinaryRelation}
Let $\preceq$ be the binary relation on $E$ where $y \preceq x$ if, for every neighborhood $U$ of $x$, there exists a neighborhood 
$V$ of $y$ and $f \in \Mod(\Sigma)$ so that $f(V) \subset U$. We say that \emph{x and y are of the same type}, denoted 
$x \sim y$, if $x \preceq y$ and $y \preceq x$, and write $E(x)$ for the set $\{ y \st y \sim x \}$.
One can easily verify that $\sim$ defines an equivalence relation.
\end{definition}

From the definition of $\preceq$, we obtain a partial order on the set of equivalence classes under $\sim$.
Indeed, the relation $\prec$, defined by $x \prec y$ if $x \preceq y$ and $x \nsim y$, gives a partial order on the set 
of equivalence classes under $\sim$.

\begin{proposition}[Proposition 4.7 in \cite{MR}]
\label{Prop:PartialOrderMaxElements}
The partial order $\prec$ has maximal elements. Furthermore, for every maximal element $x$, the equivalence class $E(x)$ 
is either finite or a Cantor set.
\end{proposition}

The above statement also applies to every clopen subset of $E$.  Let $\mathcal{M}=\mathcal{M}(E)$ denote the set of 
maximal elements for $\prec$. For a clopen subset $A \subset E$, we denote the maximal ends in $A$ by $\mathcal{M}(A)$.

\subsection{$\Map(\Sigma)$ is locally CB} A \textit{CB generating set} is the union of a CB neighborhood of the identity and a finite set. 

\begin{definition}[Locally CB]\label{Def:CoarselyBounded}
Let $G$ be a topological group. A subset $H \subseteq G$ is \textit{coarsely bounded}, abbreviated CB, in $G$ if $H$ has 
finite diameter with respect to every continuous left-invariant length function on $G$. We say $G$ is \emph{locally CB} if 
there is a neighborhood of the identity that is a CB subset of $G$. 
\end{definition}

Mann and Rafi give a classification of surfaces $\Sigma$ for which $\Map(\Sigma)$ is locally CB. 
In particular, they show that if $\Map(\Sigma)$ is locally CB, then $\Sigma$ has the following structure
\cite[Proposition 5.4]{MR}.

\begin{proposition}\label{Prop:EndSpaceDecomposition}
If $\Map(\Sigma)$ is locally CB, then there exists a subsurface of finite type $L \subset \Sigma$ giving a 
partition of the ends 
\[
E = \bigsqcup\limits_{A \in \mathcal{A}} A
\]
where $\abs{\mathcal{A}}$ equals the number of boundary components of $L$. Each set $A \in \mathcal{A}$ is clopen and 
$\mathcal{M}(A) \subset \mathcal{M}(E)$. In fact, points in $\calM(A)$ are all in the same equivalence class, and 
$\calM(A)$ is either a singleton or a Cantor set. 
\end{proposition}

Furthermore, for every $A, B \in \calA$, there is a clopen set $W_{A,B}$ such that if $E(x)$ intersects both $A$
and $B$, it also intersects $W_{A,B}$ \cite[Lemma 6.10]{MR}.

\subsection{Stable surface} In this paper, we assume the surface is stable.

\begin{definition} \label{Def:stable}
We say a neighborhood $U$ of a point $x \in E$ is \emph{stable} if  
every neighborhood $U'$ of $x$ contains a smaller neighborhood $U''$ of $x$ that is homeomorphic to $U$. The surface
$\Sigma$ is called \textit{stable} if every point $x$ in $E$ has a stable neighborhood.  
\end{definition}

In fact, a stronger conclusion holds. Namely all clopen neighborhoods of 
$x$ inside the stable neighborhood are homeomorphic.  

\begin{lemma}[Lemma 4.17 in \cite{MR}] \label{Lem:strong-stabel}
If $U$ is a stable neighborhood of $x \in E$, then for any clopen 
neighborhood $U' \subset U$ of $x$, $U'$ is homeomorphic to $U$. 
\end{lemma}

\subsection{The subgroup $\FMap(\Sigma)$}\label{SubSec:FiniteIndexSubgroup}
As mentioned above, $\calM(A)$ is either an isolated
point in $\calM(E)$ or a Cantor set. Let $\calA_{\text iso}$ be the subset of $\calA$ consisting of sets $A$ where 
$\calM(A)$ is isolated. Let $x_{A}$ denote an element of $\mathcal{M}(A)$. 
The set $A$ is always a stable neighborhood of $x_A$ and if $x_B \in E(x_A)$
then $A$ is homeomorphic to $B$. Any $f \in \Map(\Sigma)$ acts by a permutation on the set of isolated maximal ends, 
\[
M_{\text iso} = \{x_A \st A \in \calA_{\text iso} \}
\]
which is a finite set. There exists a map 
\[
\sigma \from \Map(\Sigma) \to \mathrm{Sym}(M_{\text iso})
\]
induced by the action of $\Map(\Sigma)$ on $M_{\text iso}$.
Since $M_{\text iso}$ is finite, this implies that $\mathrm{ker}(\sigma)$ 
is a finite index subgroup of $\Map(\Sigma)$. We define
$\FMap(\Sigma)=\mathrm{ker}(\sigma)$, that is, $\FMap(\Sigma)$ the subgroup 
of $\Map(\Sigma)$ that fixes the isolated maximal ends point-wise.

\subsection{Coarsely bounded generating sets}\label{SubSec:CBGeneratingSet}
In \cite{MR} Mann-Rafi gave a CB generating set for $\FMap(\Sigma)$ and $\Map(\Sigma)$. 
Let $L$ be the surface of finite-type as above. The first collection of elements in our 
generating set are the elements fixing $L$
\[
\nu_{L} = \big\{g \in \Map(\Sigma) \st g|_L = id \big\}.
\]
When $\Map(\Sigma)$ is CB generated, there is a finite set 
$F \subset \FMap(\Sigma)$ such that $\nu_L \cup F$ generates $\FMap$ 
(See \cite[Section 6.4]{MR} for details). 
We say $\pi \in \mathrm{Sym}(M_{\text iso})$ is \emph{type preserving} 
if for every $x_A \in M_{\text iso}$, $x_B=\pi(x_A) \in E(x_A)$. 
For every type preserving $\pi \in \mathrm{Sym}(M_{\text iso})$
we fix a homeomorphism $p_\pi$ which permutes the sets in $\calA$, 
sending $A$ to $B$ when $\pi(x_A) = x_B$, but fixes the finite type 
subsurface $L$ set wise (that is the restriction of $p_\pi$ to $L$ is 
an element of $\Map(L)$). In fact, we can choose these so that 
$p_{\pi^{-1}} = p_\pi^{-1}$. 
Since there are only finitely many such permutations, the collection $P = \{ p_\pi\}$ 
is finite. Then $\nu_L \cup F \cup P$ is a CB generating sets for $\Map(\Sigma)$.

\subsection{$\FMap(\Sigma)$ is quasi-isometric to $\Map(\Sigma)$.}
Let $\calS_{\FMap}$ and $\calS_{\Map}$ be the generating sets for $\FMap(\Sigma)$
and $\Map(\Sigma)$ from the previous section. Let $\norm{\param}_{\FMap}$
and $\norm{\param}_{\Map}$ be the associated word lengths for $\FMap(\Sigma)$
and $\Map(\Sigma)$, respectively. 
Since $\FMap(\Sigma)$ and $\Map(\Sigma)$ are CB generated, the word lengths
with respect to any other generating sets would be bi-Lipschitz equivalent 
to the word lengths associated to $\calS_{\FMap}$ and $\calS_{\Map}$. Hence
it is sufficient to consider these generating sets.

Also, recall that to show that a metric space $(X, d_X)$ is quasi-isometric
to its subspace $(Y, d_Y)$ it is sufficient to check that the
inclusion map $Y \to X$ is coarsely onto and a quasi-isometric embedding. 
That is there are constants $M, C>0$ such that:
\begin{itemize}
 \item for every $x \in X$ there is a $y \in Y$ such that $d_X(x,y) \leq C$. 
 \item for all $y_1, y_2 \in Y$, we have
 \[
  \frac 1M d_Y(y_1, y_2)-C \leq d_X(y_1, y_2) \leq M \cdot d_Y(y_1, y_2)-C.
 \]
\end{itemize}

\begin{theorem}\label{Thm:CompatibleGenSetQI}
Let $\Sigma$ be a surface such that $\Map(\Sigma)$ is CB generated.
Then $(\FMap, \norm{\param}_{\FMap})$ and $(\Map, \norm{\param}_{\Map})$ are quasi-isometric.
\end{theorem}

\begin{proof}
Let $\nu_L$, $F$ and $P$ be as in Section \ref{SubSec:CBGeneratingSet} so that we have 
$\calS_{\FMap} = \nu_L \cup F$ and $\calS_{\Map}= \nu_L \cup F \cup P$. 
Notice that
\[
\Map(\Sigma) = \bigcup_{p_\pi \in P} p_\pi \FMap(\Sigma).  
\]
Let $H$ be the set of elements of $\FMap(\Sigma)$ that can be written 
in the form $psp'$ where $s \in \calS_{\Map}$ and $p,p' \in P$. 
We begin by showing that
\[
\max\limits_{h \in H} \norm{h}_{\FMap}
\]
is finite. We consider the cases when $s$ is in $\nu_L$, or in $F \cup P$
separately. 

For $s \in \nu_L$ and $p, p' \in P$, if $h=p s p'$ is an element of
$\FMap(\Sigma)$, we must have $p = p_{\pi}$ and $p' = p_{\pi^{-1}}$
since $s$ acts trivially on $M_{iso}$. By the construction of $P$, 
this implies that $p' = p^{-1}$. Since, $s|_L=id$, we also have
$h|_L=id$. That is, $h \in \nu_L$ and $\norm{h}_{\FMap}=1$. 

Let $H'$ be the subset of $H$ consisting of elements of the type $h = p s p'$, 
where $s \in F \cup P$. Then $H'$ is a finite set. That is, letting 
\[
M= \max\limits_{h \in H'} \norm{h}_{\FMap}, 
\]
we also have  
\[
M = \max\limits_{h \in H } \norm{h }_{\FMap} = \max \{ M,1 \}.
\]

Now, let $g \in \FMap(\Sigma)$. We can write $g = s_1 s_2 \ldots s_n$, where 
$s_i \in \calS_{\Map}$ and $n = \norm{g}_{\Map}$. Let 
$g_i = s_is_{i+1}\dots s_n$ and let $p_i \in P$ be such that 
$g_i \in p_i \FMap(\Sigma)$, that is 
\[
g_i = s_i s_{i+1} \ldots s_n = p_{i} f_i
\]
for some $f_i \in \FMap(\Sigma)$.  Therefore, 
\[
p_{i} f_i=g_i = s_i \, g_{i+1}= s_i \, p_{{i+1}} \, f_{i+1}, 
\]
which implies
\[
f_i f_{i+1}^{-1} = p_i^{-1} s_i p_{i+1}. 
\]
Note that  $p_i^{-1},p_{i+1} \in P$, $s_i \in \calS_{\Map}$ and $f_i f_{i+1}^{-1}\in \FMap(\Sigma)$,
hence $p_i^{-1} s_i p_{i+1} \in H$ and
\[
\norm{p_i^{-1}\, s_i \, p_{i+1}}_{\FMap} \leq M. 
\]
Now, by rewriting $g$ as
\[
g = (e\, s_1\, p_2) (p_2^{-1}\,s_2\, p_3)(p_3^{-1}\,s_3\,p_4)
\ldots(p_n \,s_n\, e)
\]
we have that $p_i s_i p_{i+1} \in H$ for $1 \leq i \leq n$ (by setting $p_0=p_{n+1}=e$ where $e$
is the identity). Therefore, 
\[
\norm{g}_{\FMap} \leq M \cdot n = M \cdot \norm{g}_{\Map}. 
\]
Since $\calS_{\FMap} \subset \calS_{\Map}$, we always have $\norm{g}_{\Map}
\leq \norm{g}_{\FMap}$. That is, the word lengths 
$\norm{\param}_{\FMap}$ and $\norm{\param}_{\Map}$ are bi-Lipschitz 
equivalent in $\FMap(\Sigma)$. Also, for every $g \in \Map(\Sigma)$
there is an $f \in \FMap(\Sigma)$ and $p \in P$ such that $g=pf$, and thus
$d_{\Map}(f,g) \leq 1$. Hence $(\FMap, \norm{\param}_{\FMap})$ and 
$(\Map, \norm{\param}_{\Map})$ are quasi-isometric.
\end{proof}

\section{The Shark Tank}\label{Sec:SharkTank}
The shark tank, $T$, is a bi-infinite cylinder with a discrete countable set of punctures 
exiting in both ends. The end space $E=E(T)$ is homeomorphic to the following subset of $\RR$:
\[
\{0\} \cup \left\{ \frac 1n \right\}_{n=2}^\infty \cup \left\{ \frac {n-1}n \right\}_{n=3}^\infty \cup \{ 1 \}.  
\]
We call the ends of the cylinder which are accumulated by punctures the 
\emph{limit ends}. These are also the maximal elements of $E$. 
The group $\FMap(\Sigma)$ is the index two subgroup of $\Map(T)$ which 
fixes the limit ends point wise. In this section, we prove that $\FMap(T)$, and as a result $\Map(T)$, have infinite 
asymptotic dimension. 

The generating set given in \cite{MR} has a very simple form in the case of $\FMap(T)$ which 
we now describe. 
Fix a curve $\beta$ on $T$ that separates the two limit ends of $T$. This decomposes the end space into two sets,
which we denote by $A$ and $B$. We denote the limit end in $A$ by $x_A$ and
the limit end in $B$ by $x_B$. In addition, we fix an ordering on the punctures 
in $T$ (the non-limit ends) and label them with $p_i$, $i \in \ZZ$ such that 
$p_i \in A$ for $i \leq 0$ and $p_i \in B$ for $i >0$. 

\begin{figure}[h]
\setlength{\unitlength}{0.01\linewidth}
\begin{picture}(75,25)
\put(0,0){\includegraphics{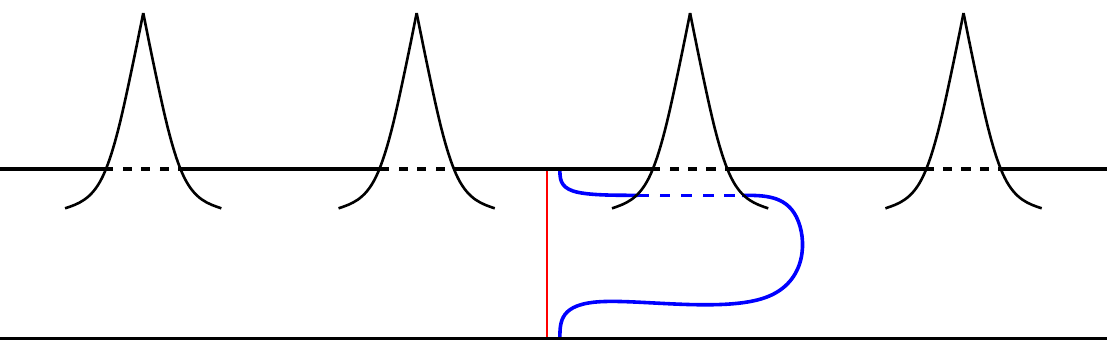}}
\put(-10,5){$x_A$}
\put(80,5){$x_B$}
\put(55,5){\textcolor{blue}{$h(\beta |_{\sigma})$}}
\put(8,24){$p_{-1}$}
\put(27,24){$p_0$}
\put(46,24){$p_1$}
\put(64,24){$p_2$}
\put(33,5){\textcolor{red}{$\beta |_{\sigma}$}}
\end{picture} 
\caption{The map $h$ is homotopic to the shift map $h_\sigma$ whose
support is the strip $\sigma$ and sends $p_i$ to $p_{i+1}$.}
\label{Fig:SharkTankBaseCurve} 
\end{figure}
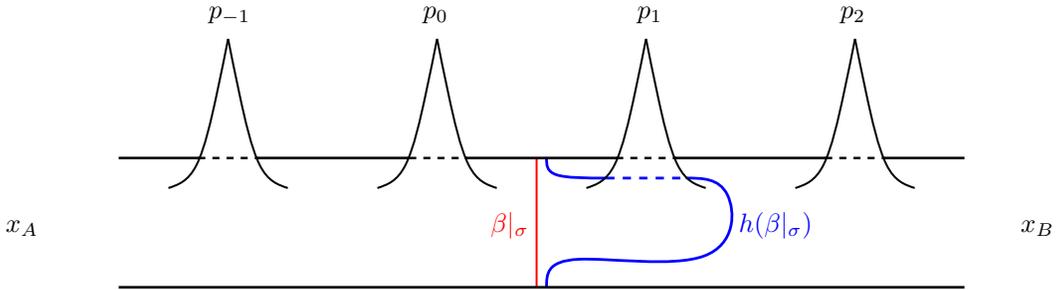

There is a homeomorphism $h \from T \to T$ fixing $x_A$ and $x_B$ such that
$h(p_i) = p_{i+1}$ and $f(\beta)$ is disjoint from $\beta$ with
$\beta$ and $f(\beta)$ bounding a a surface with two boundary components 
and one puncture.
We refer to $h$ as \emph{the shift map}. This is consistent with 
the definition of shift in the previous section; if we embed an strip $\sigma$
in $T$ containing all the punctures and limiting to $x_A$ and $x_B$, the 
homeomorphism $h_\sigma$ described in the previous section is homotopic to $h$
since the complement of $\sigma$ is a strip with no topology (see
\figref{Fig:SharkTankBaseCurve}). 

Define
\[
\nu_{\beta} = \big\{g \in \FMap(\Sigma) \st g (\beta) = \beta \big\}. 
\]
Then ${\mathscr S} =\nu_\beta \cup \{ h_\sigma \}$ is the generating set for $\FMap(T)$ given in 
\cite[Section 6.4]{MR}. We now equip $\FMap(T)$ with the word length associated to this 
generating set which we denote by $\norm{\param}_{\FMap}$. 

\begin{definition} \label{Def:length-function}
For a group $G$, a function $\norm{\param} \from G \to \mathbb{R}_{+}$ is a
called \textit{length function} if it is continuous, if $\norm{\phi}=\norm{\phi^{-1}}$ 
for all $\phi \in \FMap(\Sigma)$ and if it satisfies the triangle inequality, namely, 
for $\phi, \psi \in G$ we have
\[
\norm{\phi \, \psi} \leq \norm{\phi} + \norm{\psi}. 
\]
\end{definition}

We now define a length function on $\FMap(T)$, namely, for $\phi \in \FMap(\Sigma)$, define
\begin{equation}
\norm{\phi}= \big| \{p_i \in A \st \phi(p_i) \in B \} \big| + \big|\{p_i \in B \st \phi(p_i) \in A \} \big|. 
\end{equation}
Note that, since maps in $\FMap(T)$ fix the limit ends, they also fix 
a neighborhood of these ends and hence can move only a finitely many 
punctures from one side to another. That is, $\norm{\phi}$ is a finite number. 

\begin{theorem}\label{Thm:lfPseudoMetricSharkTank}
The function $\norm{\param}$ is a length function on $\FMap(T)$. Furthermore, 
for $\phi \in \FMap(T)$, we have 
\begin{equation} \label{Eq:Shark-norm}
\norm{\phi} \leq \norm{\phi}_{\FMap}.
\end{equation}
\end{theorem}

\begin{proof}
We start by checking the triangle inequality. 
Consider $\phi,\psi \in \FMap(T)$. For any $p_i \in A$ where $\psi \phi(p_i) \in B$,
we either have 
\[
\Big( p_i \in A \quad\text{and}\quad \quad \phi(p_i) \in B \Big) 
\qquad\text{or}\qquad
\Big( \phi(p_i) \in A \quad\text{and}\quad \psi\phi(p_i) \in B \Big).  
\]
Therefore, 
\[
\{p_i \in A \st \psi \phi (p_i) \in B \} 
\subset 
\big\{ p_i \in A \st \phi(p_i) \in B \big\} \cup  
\big\{ p_i \in A \st \phi(p_i) \in A \quad \text{and} \quad  \psi\phi(p_i) \in B \big\}.  
\]
and hence (denoting $\phi(p_i)$ with $q_i$), we have 
\[
\big| \{p_i \in A \st \psi \phi (p_i) \in B \} \big| 
\leq
\big| \big\{ p_i \in A \st \phi(p_i) \in B \big\}\big| +
\big| \big\{ q_i \in A \st \psi(q_i) \in B \big\} \big|.  
\]
Similarly, 
\[
\big| \{p_i \in B \st \psi \phi (p_i) \in A \} \big| 
\leq
\big| \big\{ p_i \in B \st \phi(p_i) \in A \big\}\big| +
\big| \big\{ q_i \in B \st \psi(q_i) \in A \big\} \big|.  
\]
Therefore, 
\[
\norm{\phi \, \psi} \leq \norm{\phi} + \norm{\psi}. 
\]
Now we check that $\norm{\param}$ is continuous. Note that, for $\phi \in \nu_\beta$, 
we have $\norm{\phi}=0$. That is $\norm{\param}$ is zero on some neighborhood
of the identity. The triangle inequality above proves that $\norm{\param}$ is
continuous. Also, since $\phi$ is a homeomorphism, 
\begin{align*}
    \norm{\phi} 
    &=\big| \{p_i \in A \st \phi(p_i) \in B \} \big| 
      + \big|\{p_i \in B \st \phi(p_i) \in A \} \big|\\
    &=\big| \{q_i \in B \st \phi^{-1}(q_i) \in A \} \big| 
      + \big|\{q_i \in A \st \phi^{-1}(q_i) \in A \} \big| = \norm{\phi^{-1}}. 
\end{align*}
Hence, $\norm{\param}$ is a length function. 

To see the second assertion, let $s \in \mathscr{S}$ be an element of the generating set
of $\FMap(T)$. If $s \in \nu_{\beta}$, then $s$ fixes the base curve $\beta$.
Therefore, $s(A)\subset A$ and $s(B) \subset B$, and therefore, $\norm{s}=0$. 
Alternatively, if $s=h_\sigma$ then $p_0$ is the only puncture in $A$ that is mapped 
to $B$ and no punctures from $B$ are mapped to $A$. That is, 
\[
\norm{h}= \big| \big\{ p_i \in A \st h(p_i) \in B \big\}\big| 
 + \big| \big\{p_i \notin B \st  h(p_i) \in B \big\}\big|= 1 + 0 =1. 
\]
Now, for $\phi \in \FMap(T)$, if $\norm{\phi}=n$, then $\phi=s_1 \dots s_n$ where 
$s_i \in\mathscr{S}$. By the triangle inequality, 
\begin{equation*} 
\norm{\phi} \leq \sum_{i=1}^n \norm{s_i} \leq n = \norm{\phi}_{\FMap}.  
\qedhere
\end{equation*}
\end{proof}

We now construct a map from $\Qinf$ to $\FMap(T)$. We begin by defining a 
homeomorphism, associated to a given element $a \in \Qinf$, that permutes 
the punctures in $B$. We make use of the following functions.
Let $z \from \Qinf \to \NN$, be the function that indicates the
number of zeros in $a$ before the final one. That is, 
\begin{equation}\label{Fn:z(a)}
z(a) = \big|\big\{ j \ \big| \ \text{$a_j=0$ 
and there is $k >j$ such that $a_k=1$} \big\} \big|.
\end{equation}
Next, we define $z_a(i)$ to be the position of $i$-th zero in the sequence 
$a$, namely, 
\begin{equation}\label{Fn:z_a}
z_{a}: \{ 1, \dots, z(a)\} \to \NN, 
\quad\text{where}\quad 
z_{a}(i)=j 
\quad\text{if $a_j=0$ \quad and }\quad
\big|\big\{ 1 \leq k \leq j \st a_k=0 \big\}\big| = j. 
\end{equation}

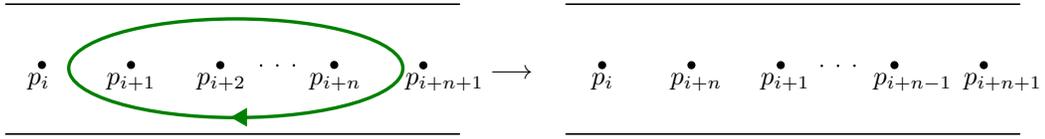
\begin{figure}[ht]
\setlength{\unitlength}{0.01\linewidth}
\begin{picture}(75,12)
\put(-5,0){\includegraphics[width=90\unitlength]{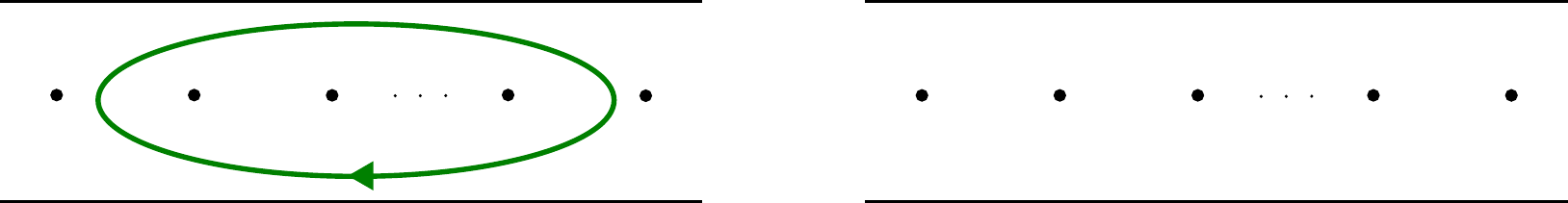}}
\put(-3,4.5){$p_i$}
\put(4,4.5){$p_{i+1}$}
\put(12,4.5){$p_{i+2}$}
\put(22,4.5){$p_{i+n}$}
\put(30.5,4.5){$p_{i+n+1}$}
\put(38,5){$\longrightarrow$}
\put(47,4.5){$p_i$}
\put(54,4.5){$p_{i+n}$}
\put(62,4.5){$p_{i+1}$}
\put(72,4.5){$p_{i+n-1}$}
\put(80,4.5){$p_{i+n+1}$}
\end{picture} 
\caption{The $1/n^{th}$-twist around the punctures $p_{i+1}, \dots, p_{i+n}$.}
\label{Fig:FractionTwist} 
\end{figure}

Now, for a given $a \in \Qinf$, we construct a \emph{puncture permutation map}
$\pi_a$ as a product of Dehn twists.  
We define a $1/n^{\text{th}}$-twist around consecutive punctures 
$p_{i+1}, \ldots, p_{i+n}$, denoted by $T^{1/n}_{{i+1}, \ldots, {i+n}}$, as follows.
Let $c$ be a curve in $\sigma$ surrounding the $n$ punctures as shown in Figure \ref{Fig:FractionTwist}. 
Then, the twist map is the homeomorphism with 
support on the surface separated by $c$ that sends $p_j$ to $p_{j+1}$ for 
$j=i+1, \dots, (n+i-1)$ and sends $p_{i+n}$ to $p_{i+1}$ such that 
$(T^{1/n}_{{i+1}, \ldots, {i+n}})^n$ is a Dehn twist around the curve $c$. 

We define the map $\pi_a$ to be 
\[
\pi_a = \prod_{i=1}^{z(a)} 
T_{{z_a(i)}, \ldots, {\abs{a}+i}}^{1/(\abs{a} + i - z_a(i))}.
\]
Observe that, the $i$--th twist moves the $(|a|+i)$-th puncture 
to the $z_a(i)$--th place. Therefore, the punctures 
$p_{\abs{a}+1}, \dots, p_{\abs{a}+ z(a)}$ are sent by $\pi_a$ to 
punctures in correspondence with the positions where $a$ has zeros.

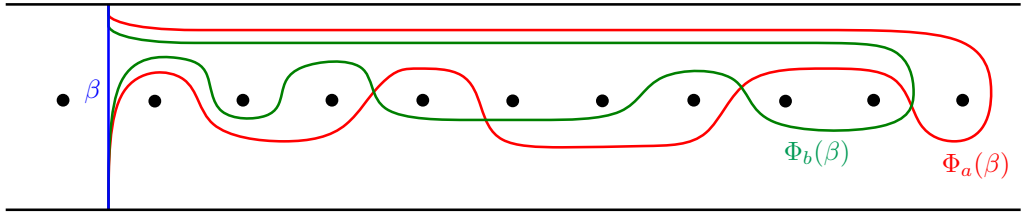
\begin{figure}[ht]
\setlength{\unitlength}{0.01\linewidth}
\begin{picture}(75,20)
\put(-7,0){\includegraphics[width=90\unitlength]{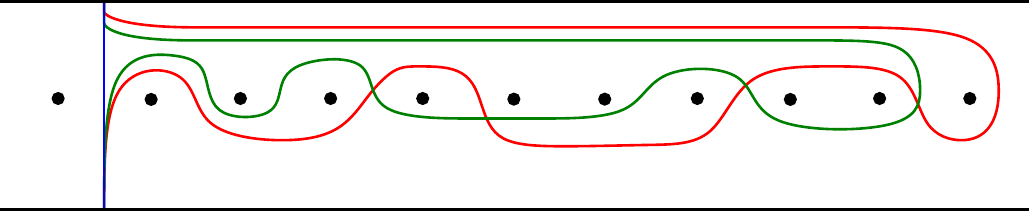}}
\put(0,10){\textcolor{blue}{$\beta$}}
\put(76,3.5){\textcolor{red}{$\Phi_a(\beta)$}}
\put(62,4.5){\textcolor{ForestGreen}{$\Phi_b(\beta)$}}
\end{picture} 
\caption{Curves $\Phi_a(\beta)$ and $\Phi_b(\beta)$ for
$a = (0,1,1,0,1,1,1,0,0,1,0,\ldots)$ and $b=(0,1,0,0,1,1,0,1,0,\ldots)$.}
\label{Fig:PuncturePermutationStripExample} 
\end{figure}

Finally, we define:
\begin{equation}\label{PhiMap}
    \begin{aligned}
    \Phi \from \Qinf &\to \FMap(T)\\
    a &\mapsto \pi_{a}h_\sigma^{\abs{a}}.
    \end{aligned}
\end{equation}
That is, we first shift $|a|$ punctures from $A$ to $B$ sending the 
punctures $p_1, \dots, p_{z(a)}$ to punctures $p_{|a|+1}, \dots, p_{|a|+z(a)}$
and then move these punctures backwards inserting them at the positions where $a$
is zero. Therefore, the punctures $p_{-|a|+1}, \dots, p_0$ are sent to the positions
where $a$ is 1. 

For any $a\in \Qinf$, the map $\pi_a h_\sigma^{\abs{a}}$ sends the 
arc $\beta|_\sigma$ to the arc $\Phi(a)(\beta|_\sigma)$ which goes around 
a puncture $p_i$ if and only if $a_i=1$, 
see Figure \ref{Fig:PuncturePermutationStripExample}.
That is, $h_\sigma^{\abs{a}}(\beta|_\sigma)$ is an arc with the same end points 
as $\beta|_\sigma$ where there are $\abs{a}$ punctures between $\beta|_\sigma$ and 
$h^{\abs{a}}(\beta|_\sigma)$, which are exactly the punctures $p_i$ where $a_i=1$.

\begin{figure}[ht]
\setlength{\unitlength}{0.01\linewidth}
\begin{picture}(75,27)
\put(3,0){\includegraphics[width=70\unitlength]{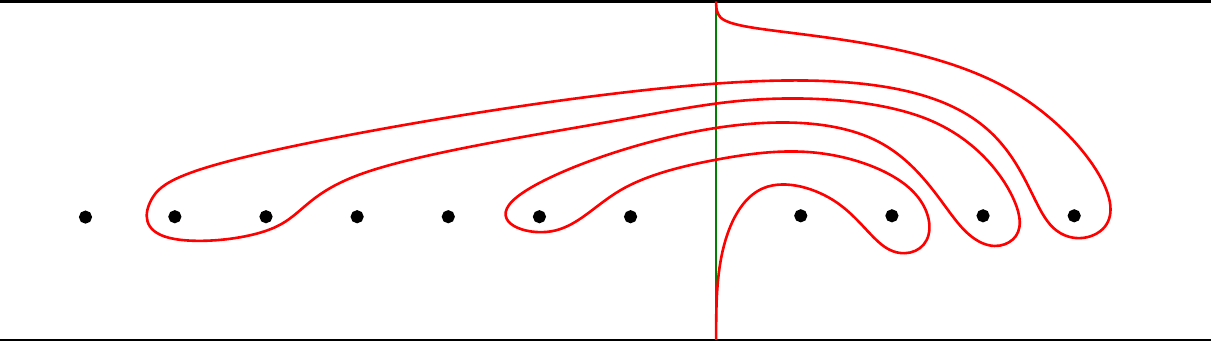}}
\put(69,10){$\Phi(a) \Phi^{-1}(b)(\beta|_\sigma)$}
\end{picture} 
\caption{There are $|a-b|$ punctures between $\beta|_\sigma$ and
$\Phi(a)\Phi_b^{-1}(\beta|_\sigma)$ one for every index $i$ where $a_i \not = b_i$.}
\label{Fig:PuncturePermutationStripExampleTaut} 
\end{figure}

For $a, b \in \Qinf$, we would like to compute $\norm{\Phi(b)^{-1}\Phi(a)}$. 
For any $i>0$ where $a_i=1$, the map $\Phi(a)$ moves a puncture from 
$A$ into $B$ placing it at $p_i$. Similarly, for any $i>0$ where $b_i=1$, 
the map $\Phi(b)$ moves a puncture from $A$ to $B$ placing it at $p_i$.
If $a_i=b_i=1$ the puncture from $A$ that is sent to $p_i$ by $\Phi(a)$ 
is sent back to $A$ by $\Phi(b)^{-1}$. 
If $a_i=1$ and $b_i=0$, a puncture from $A$ is sent to $p_i$ and remains 
in $B$ after applying $\Phi(b)^{-1}$. If $a_i=0$ and $b_i=1$, the
the puncture that was sent to $p_i$ by $\Phi(a)$ was in B, but then 
$\Phi(b)^{-1}$ sends $p_i$ to $A$. 

To visualize this, consider two arcs $\Phi(a)(\beta|_\sigma)$ and $\Phi(b)(\beta|_\sigma)$ 
on the strip $\sigma$ (see Figure \ref{Fig:PuncturePermutationStripExample}). 
Applying $\Phi(b)^{-1}$ to the arc $\Phi(a)(\beta|_\sigma)$ is like 
pulling the arc $\Phi(b)(\beta)$ taut so that it is back in alignment with the arc $\beta$. 
Then every puncture that was between $\Phi(b)(\beta|_\sigma)$ and $\beta|_\sigma$ will end
up in $A$ (See Figure \ref{Fig:PuncturePermutationStripExampleTaut}).
There are two sets of punctures between $\beta$ and 
$\Phi(b)^{-1}\Phi(a)(\beta)$; one puncture associated to every $i$ where 
$a_i=0$ and $b_i=1$ to the left of $\beta$, and one puncture for every 
$i$ where $a_i =1$ and $b_i=0$ to the right of $\beta$. That is,
a puncture is moved from $A$ to $B$ or from $B$ to $A$ for every $i>0$ where
$a_i \not = b_i$. This implies that 
\begin{equation}\label{Eqn:QinfNormRelPhi}
\norm{\Phi(b)^{-1}\Phi(a)(\beta)} = |a-b|. 
\end{equation} 

We now show that $\Phi$ a quasi-isometric embedding from $\Qinf$ to $\FMap (T)$. 
Recall that the distance between two elements $f, g, \in \FMap(T)$ is
$\norm{fg^{-1}}_{\FMap}$, where $\norm{\param}_{\FMap}$ is the word
length with respect to the generating set $\calS$. A quasi-isometric embedding
is a map that preserve distances up to uniform additive and multiplicative errors. 

\begin{proposition}\label{Prop:QIsometryQinfPhiSharkTank}
For $\Phi$ as in Equation \ref{PhiMap} and $a,b \in \Qinf$,  we have  
\[
\abs{a-b} \leq \norm{\Phi(b)^{-1}\Phi(a)}_{\FMap} \leq \abs{a-b}+3.
\]
That is, the map 
\[
\Phi \from (\Qinf, \abs{\param}) \to (\FMap ( T), \norm{\param}_{\FMap})
\]
is a quasi-isometric embedding. 
\end{proposition}

\begin{proof}
The left inequality follows from Equations \eqref{Eq:Shark-norm} and
\eqref{Eqn:QinfNormRelPhi}. To prove the right  inequality, 
let $a, b \in \Qinf$, we will find elements in $s_i \in \calS$ 
such that $(\prod s_i )\circ \Phi(b)^{-1}\Phi(a)$ is the identity.

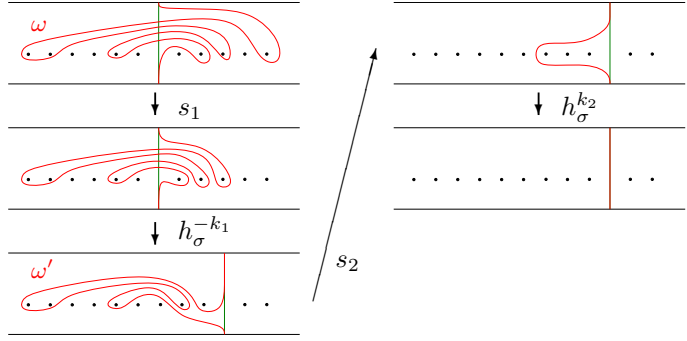
\begin{figure}[ht]
\setlength{\unitlength}{0.01\linewidth}
\begin{picture}(60,32)
\includegraphics[width=0.6\textwidth]{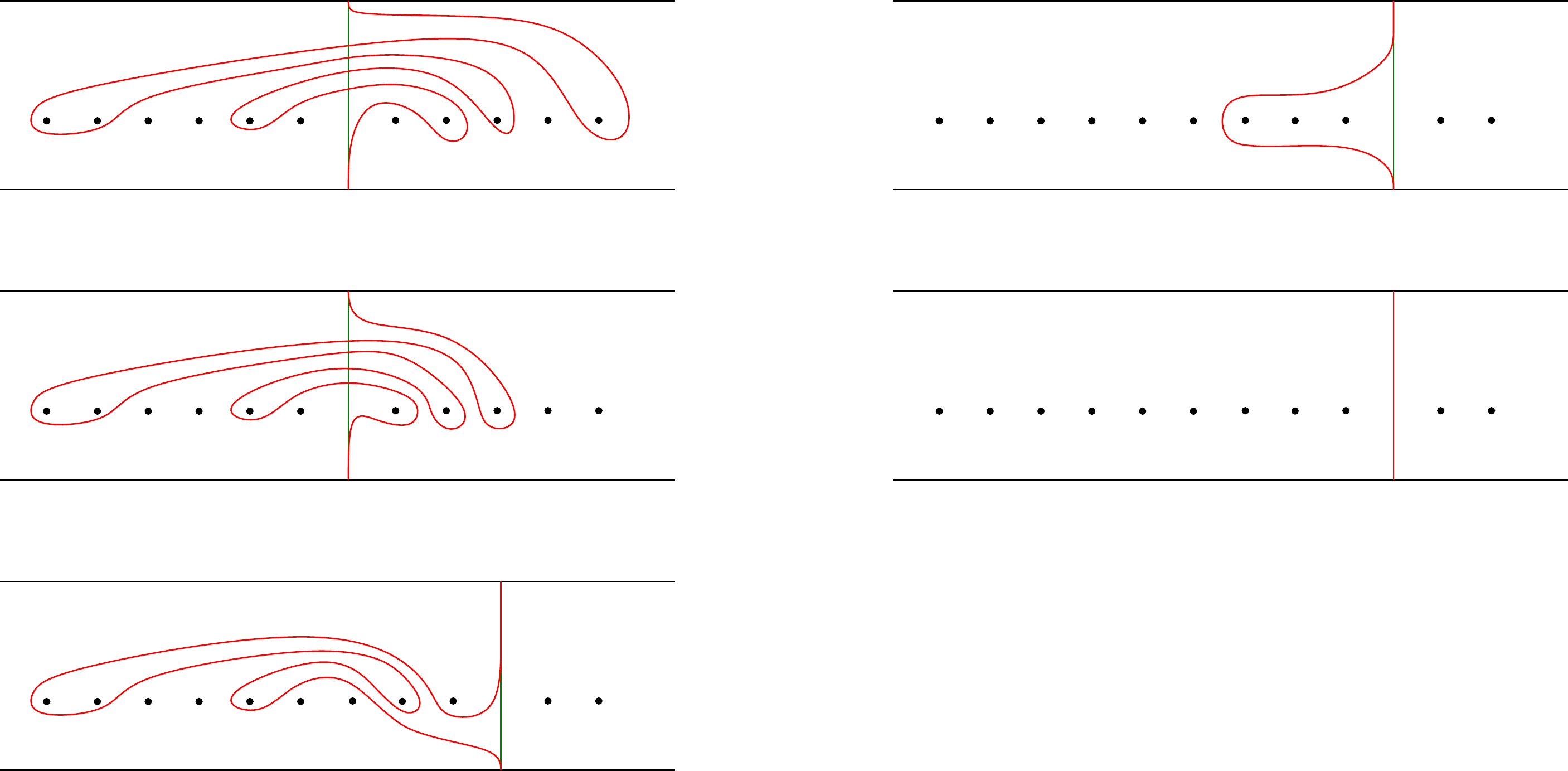}
\linethickness{.7pt}
\put(-47,21.5){\vector(0,-1){2}}
\put(-45,19.5){$s_1$}
\put(-47,10){\vector(0,-1){2}}
\put(-45,8.5){$h_\sigma^{-k_1}$}
\put(-13,21.5){\vector(0,-1){2}}
\put(-11,19.5){$h_\sigma^{k_2}$}
\put(-33, 3){\vector(1,4){5.6}}
\put(-31,6){$s_2$}
\put(-58,27){\textcolor{red}{\small{$\omega$}}}
\put(-58,5){\textcolor{red}{\small{$\omega'$}}}
\end{picture}
\caption{The arc $\Phi(b)^{-1}\Phi(a)(\beta|_\sigma)$ can be sent to 
$\beta|_\sigma$ in five steps.}
\label{fig:PuncturePermutationInverse}
\end{figure}

Let $\omega= \Phi(b)^{-1}\Phi(a)(\beta|_\sigma)$. There are
$k_1$ punctures between $\omega$ and $\beta|_\sigma$ in $B$, one for every index $i$ where 
$a_i=1$ and $b_i=0$, and $k_2$ punctures between $\omega$ and
$\beta|_\sigma$ in $A$, one for every index $i$ where $b_i=1$ and $a_i=0$. 
Using an element $s_1 \in \nu_\beta$ we can 
line up the $k_1$ punctures to be in the positions of 
$p_1, p_2, ... p_{k_2}$. Then $h_\sigma^{-k_1}$ sends all these punctures
to $A$. That is, $h_\sigma^{-k_1} s_1$ sends $\omega$ to an arc $\omega'$
that is completely to the left side of $\beta$. There are still 
$k_2$ punctures between $\omega'$ and $\beta|_\sigma$. We now find
$s_2 \in \nu_\beta$ that lines up these punctures to the position 
$p_{-k_2+1}, \dots, p_0$. Then $h_\sigma^{k_2}$ sends these punctures
back to $B$. That is $h_\sigma^{k_2} s_2(\omega')= \beta|_\sigma$. 
Now the composition 
\[
\big( h_\sigma^{k_2} s_2 h_\sigma^{-k_1} s_1 \big) \circ \Phi(b)^{-1}\Phi(a)
\]
sends every puncture in $A$ to a puncture in $A$ and every puncture in $B$ to a puncture in $B$. 
Hence, this map is equal to some element $s_3^{-1} \in \nu_\beta$. Therefore, we have shown that
\[
\big(s_3 h_\sigma^{k_2} s_2 h_\sigma^{-k_1} s_1\big)\circ \Phi(b)^{-1}\Phi(a) = id. 
\]
Hence 
\[
\norm{\Phi(b)^{-1}\Phi(a)}_\calS \leq k_1 + k_2 + 3,
\]
where $k_1 + k_2 = |a-b|$.
This completes the proof. 
\end{proof}

\begin{theorem}\label{Thm:FISharkTankInfiniteAsympDim}
The group $\FMap(T)$ has infinite asymptotic dimension.
\end{theorem}

\begin{proof}
By Proposition, \ref{Prop:QIsometryQinfPhiSharkTank} the map $\Phi \from \Qinf \to \FMap(T)$ is 
a quasi-isometric embedding. Since $\Qinf$ has infinite asymptotic dimension 
\thmref{Thm:Infinite-Asymptotic-Dimension-Qinf}, implies that $\FMap(T)$ also has infinite asymptotic dimension.
\end{proof}

By Theorem \ref{Thm:CompatibleGenSetQI}, $(\FMap(T),\norm{\param}_\calS)$ and $(\Map(T),\norm{\param}_{\calS_{\Map}})$ are quasi-isometric.
Fact \ref{Fact:EmbeddingAsympDim} tells us that asymptotic dimension is preserved under quasi-isometry.
Therefore since $(\FMap(T),\norm{\param}_\calS)$ has infinite asymptotic dimension, so does $(\Map(T),\norm{\param}_{\calS_{\Map}})$.

\begin{theorem}\label{Thm:SharkTankInfiniteAsympDim}
The group $\Map(T)$ has infinite asymptotic dimension.
\end{theorem}

\section{The length functions} \label{Sec:Length-Function}

As we saw in the last section, the length function $\norm{\param}$, defined on $\FMap(T)$, 
is how we provide a lower bound for the word length. In this section, we show that when 
the end space is two-sided, there is a similar length function on $\FMap(\Sigma)$ that is 
bounded above by the word length. This implies that the word lengths of powers of an associated shift map 
grow linearly, and hence the shift map is essential. 

\subsection{Two-sided end}\label{SubSec:TwoSidedEnd}
Assume, for $z \in E$, that $E(z)$ is two-sided. Recall from the introduction that this means 
$\Accu(z) = X \sqcup Y$, where $X, Y$ are non-empty disjoint closed $\FMap(\Sigma)$--invariant subsets of $E$ and
where $\Accu(z)$ is the accumulation set of $E(z)$.  

Let the subsurface $L$ be as in \propref{Prop:EndSpaceDecomposition}
giving the decomposition 
\[
E = \bigsqcup\limits_{A \in \mathcal{A}} A. 
\]
As before, we fix a point $x_A \in \calM(A) \subset \calM(E)$. 

Observe that, for $A \in \calA$, either $X \cap A = \emptyset$ or $Y \cap A = \emptyset$. 
This is because, if $x \in X \cap A$, then $E(x) \subset X$ and $\Accu(x) \subset X$. But $x_A$ is 
an accumulation point of every type of point in $A$. That is, $x_A \in \Accu(x) \subset X$. 
Similarly, if $y \in Y \cap A$, then $x_A \in \Accu(y) \subset Y$. This contradicts the assumption that 
$X$ and $Y$ are disjoint. Hence, we can find a curve $\beta$ in $L$ so that $\Sigma \setminus \beta$ consists of 
two subsurfaces $\Sigma_{+}$ and $\Sigma_{-}$ such that $X \subset E(\Sigma_{-})$ and $Y \subset E(\Sigma_{+})$.
Then there is a decomposition $\calA = \calA_+ \sqcup \calA_-$ such that 
$E(\Sigma_+) = \cup_{A \in \calA_+} A$ and $E(\Sigma_-) = \cup_{B \in \calA_-} B$. 
We denote $E(\Sigma_+)$ by $E_+$ and $E(\Sigma_-)$ by $E_-$. 

\begin{theorem}[Length function associated to a two-sided end]\label{Thm:z-length-function}
For $z \in E$, assume $E(z)$ is two-sided. Let $\beta$ be the curve defined above so that 
\[
\Sigma \setminus \beta = \Sigma_{-} \cup \Sigma_{+}, \qquad
X \subset E_-=E(\Sigma_{-})\qquad\text{and}\qquad
Y \subset E_+=E(\Sigma_{+}).
\] 
Then, the function $\norm{\param} \from \FMap(\Sigma) \to \mathbb{Z}$ defined by
\begin{equation} \label{Eq:end-norm}
    \norm{\phi} = \big| \{p \in E(z) \st p \in E_{-}, \quad \phi(p) \in E_{+} \} \big| 
    + \big| \{p \in E(z) \st p \in E_{+}, \quad \phi(p) \in E_{-} \} \big|
\end{equation}
is a length function on $\FMap(\Sigma)$.
\end{theorem}

\begin{proof} 
We  first check that, for every $\phi \in \FMap(\Sigma)$, $\norm{\phi}$ is finite. 
Suppose towards a contradiction that $\norm{\phi}$ is infinite. By replacing $\phi$ with $\phi^{-1}$
if necessary, we may assume that $\phi$ maps infinitely many points from $E(z) \cap E_{-}$ into 
$E_{+}$. That is, there is a sequence $z_i \in E(z) \cap E_{-}$ such that $\phi(z_i)\in E_{+}$. 
After taking a sub-sequence, we can assume $z_i \to x \in X$. Since $\phi$ is continuous, 
and $E_{+}$ is closed, $\phi(x) \in E_{+}$. But this contradicts the fact that $X \subset E_-$ is $\FMap(\Sigma)$ 
invariant.

The proof of the triangle inequality and the fact that $\norm{\phi} = \norm{\phi^{-1}}$
are identical to the proof in Theorem \ref{Thm:lfPseudoMetricSharkTank}. 
Also, the subgroup of $\FMap(\Sigma)$ that fixes $L$ is a neighborhood of
the identity where $\norm{\param}$ is zero, hence $\norm{\param}$ is continuous. 
\end{proof}

Proceeding as in Section \ref{Sec:SharkTank}, we now show that the length function $\norm{\param}$ is
bounded above by a uniform multiple of the word length. 

\begin{theorem}\label{Thm:LRelationToWordMetricEnds}
For $z \in E$, assume $E(z)$ is two-sided and let $\norm{\param}$ be the associated length function. Let $\calS$ be a CB generating set for
$\FMap(\Sigma)$, and let $\norm{\param}_{\calS}$ denote the associated word length on $\FMap(\Sigma)$. 
Then there exists a constant $c>0$ such that for every $\psi \in \FMap(\Sigma)$
\begin{equation}
\norm{\psi} \leq c \cdot \norm{\psi}_{\calS}.
\end{equation}
\end{theorem}

\begin{proof}
Since $\FMap(\Sigma)$ is $CB$ generated, the word length associated to every two CB generating sets are Lipschitz 
equivalent. Hence, without loss of generality, we can assume $\FMap(\Sigma)$ is equipped with the generating set given 
in Section \ref{SubSec:CBGeneratingSet}. That is, there is a finite set $F$ such that $\calS = \nu_L \cup F$. 
For $A \in \calA$, denote the complementary component of $(\Sigma - L)$ associated to $A$ by $\Sigma_A$. 
We further write 
\[
\nu_L = \cup_{A \in \calA} \nu_A
\]
where elements of $\nu_A$ have support in $\Sigma_A$.

If $s \in \nu_L$, then $s$ fixes every $A \in \calA$ set-wise. In particular, $s$ fixes $E_+$ and $E_-$
set-wise. Therefore, $\norm{s}=0$. Now define 
\[
c = \max_{s \in F} \norm{s}. 
\]
Then, for $s \in \calS$, $\norm{s}\leq c$. For $\phi \in \FMap(\Sigma)$, if 
$\phi = s_1 \dots s_n$, 
where $s_i \in S$ and $n = \norm{\phi}_{\calS}$, then 
\[
\norm{\phi} \leq \sum_{i=1}^n \norm{s_i} \leq c \cdot n.
\] 
That is, $\norm{\phi} \leq c \cdot  \norm{\phi}_{\calS}$. 
\end{proof}

\subsection{Two-sided $E^G$}
Assume $E^G$ is two-sided. Recall from the introduction that this implies that $E^G = X \sqcup Y$, where $X, Y$ are 
non-empty disjoint closed $\FMap(\Sigma)$--invariant subsets of $E$. As in the previous section, for every 
$A \in \calA$, we have either $X \cap A = \emptyset$ or
$Y\cap A =\emptyset$. Therefore, we can choose a curve $\beta$ in $L$ giving a decomposition  
$\Sigma = \Sigma_- \cup \Sigma_+$, where $X \subset E_- = E(\Sigma_-)$ and $Y \subset E_+ = E(\Sigma_+)$.

Heuristically, to follow the two-sided end case, for $\phi \in \FMap(\Sigma)$ we would like to count the genus of the 
subsurface of $\Sigma_-$ that is moved by $\phi$ to $\Sigma_+$ plus the genus of the subsurface of $\Sigma_+$ 
that is moved by $\phi$ to $\Sigma_-$. But this is not the correct measurement. For example, consider a large 
genus subsurface $Y$ that intersects both $\Sigma_-$ and $\Sigma_+$ and let $\phi$ be a pseudo-Anosov 
homeomorphism with support in $Y$. Then no subsurface of $\Sigma_-$ is moved to $\Sigma_+$. 
Instead, $\norm{\phi}$ will be the genus of the subsurface $Y$. 

We use the $\ZZ_2$--homology of the surface $\Sigma$, where $\ZZ_2=\ZZ/2\ZZ$. 
We note that separating curves on infinite-type surfaces represent elements of homology which are not necessarily 
zero (see \cite{PV}). However, we quotient $H_1(\Sigma, \ZZ_2)$ by $H_{1}^{sep}(\Sigma, \ZZ_2)$, the subgroup 
of $H_1(\Sigma, \ZZ_2)$ generated by homology classes that can be represented
by simple separating closed curves on the surface. We denote this quotient by 
\[
\widehat{H} = \widehat{H}_1(\Sigma, \ZZ_2) = H_1(\Sigma, \ZZ_2)/H_{1}^{sep}(\Sigma, \ZZ_2).
\]
Since the curve $\beta$ represents a trivial class in $\widehat H$, the decomposition $\Sigma = \Sigma_+ \cup \Sigma_-$ 
gives a splitting $\widehat H = H_+ \oplus H_-$ where 
\[
H_{+} = \widehat{H}_1(\Sigma_{+}, \mathbb{Z}/2)
\qquad\text{and}\qquad 
H_{-} = \widehat{H}_1(\Sigma_{-}, \mathbb{Z}/2).
\]
Note that $H_- \cap H_+ = \emptyset$ since $\Sigma_-$ and $\Sigma_+$ are disjoint.

For $\phi \in \FMap(\Sigma)$, we denote the induced map on homology by 
$\phi^{*} \from H_1(\Sigma, \mathbb{Z}/2) \to H_1(\Sigma, \mathbb{Z}/2)$. Since $\phi$ sends separating 
curves to separating curves, we also have an induced isomorphism $\phi^* \from \widehat H \to \widehat H$. 

Define $\norm{\param} \from \FMap(\Sigma) \to \mathbb{Z}_{+}$ by 
\begin{equation} \label{Eq:homology-norm}
\norm{\phi} = \codim \Big( \widehat H : \big(H_+ \cap \phi^*(H_+)\big) \oplus \big(H_- \cap \phi^*(H_-)\big) \Big),
\end{equation}
where $\codim(V : W)$ is the co-dimension of a subspace $W$ in a vector space $V$.  

\begin{theorem}\label{Thm:E^G-length-function}
Assume $E^G$ is two-sided. Then $\norm{\param} \from \FMap(\Sigma) \to \NN$ is a length function. 
\end{theorem}

\begin{proof}
We start by showing that, for $\phi \in \FMap(\Sigma)$, $\norm{\phi}$ is finite. 
Since $\phi$ fixes $X$ set-wise, there is a small neighborhood of $X$ in $\Sigma_-$ that is contained in 
$\Sigma_-$. That is, there is a subsurface $\Sigma' \subset \Sigma_-$ such that $X \subset E(\Sigma')$ and 
$\phi(\Sigma') \subset \Sigma_-$. The subsurfaces $(\Sigma_- - \Sigma')$ and $(\Sigma_- - \phi(\Sigma'))$ both have 
finite genus because the fact that $E^G \cap E(\Sigma_-) = X$
implies that their end space is disjoint from $E^G$ 
This means that 
\[
\codim (H_-: \widehat H_1(\Sigma', \ZZ_2)) < \infty
\qquad\text{and}\qquad 
\codim (H_-:  \widehat H_1(\phi(\Sigma'), \ZZ_2) ) < \infty. 
\]
But
\[
\widehat H_1(\Sigma', \ZZ_2) \cap \widehat H_1(\phi(\Sigma'), \ZZ_2)  \subset H_- \cap \phi^*(H_-). 
\]
Therefore, 
\[
\codim (H_-: H_- \cap \phi^*(H_-)) \leq 
\codim (H_-: \widehat H_1(\Sigma', \ZZ_2)) + 
\codim (H_-:  \widehat H_1(\phi(\Sigma'), \ZZ_2) ) < \infty. 
\]
Similarly, $\codim (H_+: H_+ \cap \phi^*(H_+)) < \infty$ and hence $\norm{\phi} < \infty$. 

We now check that $\norm{\param}$ satisfies the triangle inequality. Consider $\phi, \psi \in \FMap(\Sigma)$. 
We have 
\begin{align*}
\codim\big(H_\pm &: H_\pm \cap \psi(\phi(H_\pm)) \big)  
\leq \codim\big(H_\pm: H_\pm \cap \psi(H_\pm) \cap \psi(\phi(H_\pm)) \big) \\
&= \codim\big(H_\pm: H_\pm \cap \psi(H_\pm) \big) 
  + \codim\big(H_\pm \cap \psi(H_\pm) : H_\pm \cap \psi(H_\pm) \cap \psi(\phi(H_\pm)) \big)\\
&\leq \codim\big(H_\pm: H_\pm \cap \psi(H_\pm) \big) 
  + \codim\big(\psi(H_\pm) : \psi(H_\pm) \cap \psi(\phi(H_\pm)) \big)\\
&= \codim\big(H_\pm: H_\pm \cap \phi(H_\pm) \big) 
  + \codim\big(H_\pm : H_\pm \cap \phi(H_\pm) \big).
\end{align*} 
Summing over $\pm$ we get 
\[
\norm{\psi \circ \phi} \leq \norm{\psi} + \norm{\phi}. 
\]

The subgroup $\nu_L$ of $\FMap(\Sigma)$ is a neighborhood of the identity and for $\phi \in \nu_L$,
$\phi$ preserves $\Sigma_\pm$. Hence, $\phi^*$ fixes $H_\pm$ and $\norm{\phi}=0$.
This fact and the triangle inequality imply that $\norm{\param}$ is continuous. Also, since 
\[
(\phi^*)^{-1} (H_\pm \cap \phi^*(H_\pm)) = (\phi^*)^{-1} (H_\pm) \cap (H_\pm)
\]
we have $\norm{\phi} = \norm{\phi^{-1}}$. Therefore $\norm{\param}$
is a length function. 
\end{proof}

We now show that the function $\norm{\param}$ is bounded above by the word metric on $\FMap(\Sigma)$.

\begin{theorem}\label{Thm:HandleShiftBoundedNorm}
There exist a constant $c >0$ such that 
\begin{equation}
\norm{\phi} \leq c \cdot \norm{\phi}_{\calS},
\end{equation}
where $\norm{\param}_\calS$ denotes the word metric on $\FMap(\Lambda)$.
\end{theorem}

\begin{proof}
This proof follows the same reasoning as the proof of Theorem \ref{Thm:LRelationToWordMetricEnds}, 
except we use that $E^G$ is 2-sided in place of $E(z)$ being two-sided.
As before, we assume that $\FMap(\Sigma)$ is equipped with the generating set $(\nu_L, F)$ given in Section 
\ref{SubSec:CBGeneratingSet}. If $s\in \nu_L$, then $s$ fixes 
the curve $\beta$ which is contained in $L$ and hence it fixes $\Sigma_\pm$. 
In particular, $s^*$ fixes $H_\pm$. Therefore, $\norm{\phi}=0$. Now, for 
\[
c = \max_{\phi \in F} \norm{\phi},
\]
we have $\norm{\phi} \leq c \cdot \norm{\phi}_{\calS}$. 
\end{proof}

\section{Infinite asymptotic dimension}\label{Sec:InfAsympDimGeneral}

In this section, we prove that the mapping class group of an infinite type surface $\Sigma$ with a two-sided end space 
has infinite asymptotic dimension.  For any such surface, we construct an associated shift map and
we use the length function from \secref{Sec:Length-Function} to show that they are essential. We then follow 
the arguments in \secref{Sec:SharkTank} to embed $\Qinf$ in $\FMap(\Sigma)$.

\begin{proposition}  \label{Prop:2-sided-implies-essential}
If $E=E(\Sigma)$ is two-sided, then $\FMap(\Sigma)$ contains an essential shift. 
\end{proposition} 

\begin{proof}
Assume first that $E^G$ is two-sided. That is, $E^G = X \sqcup Y$ where $X$ and $Y$ are non-empty closed
$\FMap(\Sigma)$ invariant sets. Let $\beta$ be a curve in $L$ separating $X$ from $Y$ and, as before, 
let $\Sigma_-$ and $\Sigma_+$ be the components of $\Sigma-\beta$. Also, let $\norm{\param}$ be the norm 
defined in \thmref{Thm:E^G-length-function}. 

Pick $x \in X$ and $y \in Y$. Since every neighborhood of $x$ and $y$ in $\Sigma$
have non-zero genus, we can find a sequence of disjoint surfaces $\Sigma_i$, each homeomorphic to a genus
one surface with one boundary component, such that $\Sigma_i \to x$ as $i \to -\infty$ and $\Sigma_i \to y$
as $i \to \infty$. We choose these so that $\Sigma_i \in \Sigma_-$ for $i \leq 0$ and $\Sigma_i \in \Sigma_+$ 
for $i > 0$. 

We connect $\Sigma_{i-1}$ to $\Sigma_i$ by an arc $\omega_i$ so that the arcs $\omega_i$ are
disjoint from each other and the other $\Sigma_j$, and so only $\omega_0$ intersects $\beta$. 
Let $\sigma$ be a regular neighborhood of the union of the $\Sigma_i$ and $\omega_i$. Then $\sigma$ is a strip 
of infinite genus exiting towards $x$ and $y$ in each end and intersecting $\beta$ in an arc $\beta|_\sigma$. 
Denote the components of $\sigma- \beta|_\sigma$ by $\sigma_-$ and $\sigma_+$. 

Let $h_\sigma$ be the shift map with support in $\sigma$ sending $\Sigma_i$ to $\Sigma_{i+1}$. 
Choosing a basis for the homology of each $\Sigma_i$ and extending it to the homology of $\Sigma_-$
and $\Sigma_+$, we can decompose the homology of $\Sigma_\pm$ as follow:
\[
H_- = \overline H_- \oplus \bigoplus_{i=-\infty}^{0} H_i
\qquad\text{and}\qquad
H_+ = \overline H_+ \oplus \bigoplus_{i=1}^{\infty} H_i
\]
where $H_i = \widehat H(\Sigma_i, \ZZ_2)$ and $\overline H_\pm$ are the homology of $\Sigma_\pm -\sigma$. 
Then, $h_\sigma^*$ fixes $\overline H_\pm$ and sends $H_i$ to $H_{i+1}$. Therefore, 
for $n>0$
\[
 (h_\sigma^*)^n(H_-) = H_- \oplus \bigoplus_{i=1}^{n} H_i.
\qquad\text{and}\qquad
(h_\sigma^*)^n(H_+) = \overline H_+ \oplus \bigoplus_{i=n+1}^{\infty} H_i
\]
That is, 
\[
\norm{h_\sigma^n} = 
\codim \Big( \widehat H : \big(H_+ \cap (h_{\sigma}^*)^n(H_+)\big) \oplus \big(H_- \cap (h_{\sigma}^*)^n(H_-)\big) \Big)
= \dim\left( {\bigoplus_{i=1}^{n} H_i}\right) =2n. 
\]
\thmref{Thm:HandleShiftBoundedNorm} implies that the diameter of the group $\langle h_\sigma \rangle$ is infinite. 
Hence, by definition, $\langle h_\sigma \rangle$ is not CB and hence $h_\sigma$ is essential. 

Now assume $E(z)$ is two-sided. That is, $\Accu(z)= X \sqcup Y$ where $X$ and $Y$ are non-empty closed
$\FMap(\Sigma)$ invariant sets. The proof in this case is similar. Let $\beta$ be a curve in $L$ separating $X$ from 
$Y$, let $\Sigma_-$ and $\Sigma_+$ be the components of $\Sigma-\beta$, and let $\norm{\param}$ be the norm 
defined in \thmref{Thm:z-length-function}. Pick $x \in X$ and $y \in Y$. Since every neighborhood of $x$ and $y$ 
in $\Sigma$ has a point in $E(z)$, we can find a sequence of points $z_i \in E(z)$ where $z_i \to x$
as $i \to -\infty$ and $z_i \to y$ as $i \to \infty$, where $z_i \in E_-$ for $i \leq 0$ and $z_i \in \Sigma_+$ 
for $i > 0$. Choose disjoint stable neighborhoods $E_i$ of $z_i$ so that $\calM(E_i) = z_i$; this is possible
since $E(z)$ is countable. Let $\Sigma_i$ be a subsurface of $\Sigma$ with one boundary component and with  $E(\Sigma_i)=E_i$. We can assume the $\Sigma_i$ are disjoint from each other and from $\beta$. 
Then $\Sigma_i \to x$ as $i \to -\infty$, $\Sigma_i \to y$ as $i \to \infty$, $\Sigma_i \in \Sigma_-$ for $i < 0$, 
and $\Sigma_i \in \Sigma_+$ for $i \geq 0$. 

As before, connect $\Sigma_{i-1}$ to $\Sigma_i$ by an arc $\omega_i$ so that the $\omega_i$ are
disjoint from each other and the other $\Sigma_j$, and so only $\omega_0$ intersects $\beta$. 
Let $\sigma$ be a regular neighborhood of the union of $\Sigma_i$ and $\omega_i$. Then $\sigma$ is a strip 
containing $\Sigma_i$ exiting towards $x$ and $y$ in each end and intersecting $\beta$ in an arc $\beta|_\sigma$. 
Then, $h_\sigma$ acts as a shift on $z_i$ with $h_\sigma^n(z_i)=z_{i+n}$. Therefore, 
\[
\norm{h_\sigma^n}=n, 
\]
which implies  $\langle h_\sigma \rangle$ has infinite diameter and $h_\sigma$ is essential. 
\end{proof}

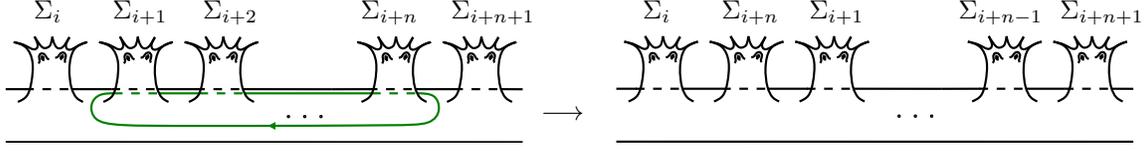
\begin{figure}[ht]
\setlength{\unitlength}{0.01\linewidth}
\begin{picture}(75,13)
\put(-12.5,0){\includegraphics[width=100\unitlength]{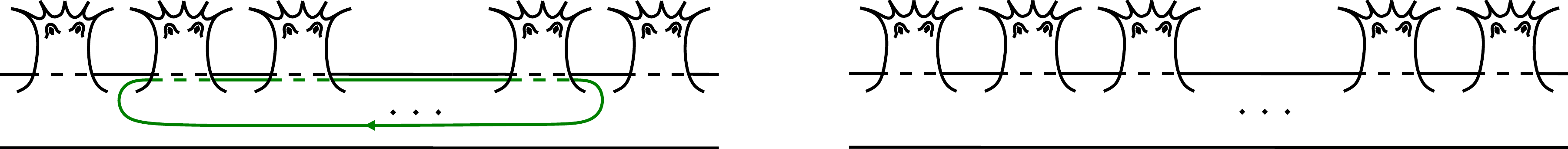}}
\put(-10,11){$\Sigma_i$}
\put(-3,11){$\Sigma_{i+1}$}
\put(5,11){$\Sigma_{i+2}$}
\put(19,11){$\Sigma_{i+n}$}
\put(27,11){$\Sigma_{i+n+1}$}
\put(35,2){$\longrightarrow$}
\put(44,11){$\Sigma_i$}
\put(51,11){$\Sigma_{i+n}$}
\put(58.75,11){$\Sigma_{i+1}$}
\put(72,11){$\Sigma_{i+n-1}$}
\put(81,11){$\Sigma_{i+n+1}$}
\end{picture} 
\caption{A $1/n^{th}$-twist around surfaces $\Sigma_i, \dots, \Sigma_{i+n}$.}
\label{Fig:FractionTwistSubSurfaces} 
\end{figure}

Similar to the construction of a puncture permutation map from Section \ref{Sec:SharkTank}, we now define a 
\textit{subsurface permutation map} for the strip $\sigma$ constructed above
where we replace the punctures with subsurfaces $\Sigma_i$. 
Again, we make use of the function $z \from \Qinf \to \ZZ$ defined in \eqnref{Fn:z(a)} which, for $a \in \Qinf$, gives 
the number of zeros in $a$ before the final one and the function $z_{a}\from \{1,...,z(a) \} \rightarrow \NN$ 
defined in the equation \eqnref{Fn:z_a} which gives the positions of the zeros in $a$. The homeomorphism $\pi_a$ is a 
composition of twists on the infinite strip $\sigma$. We define a $1/n^{\text{th}}$-twist around consecutive subsurfaces 
$\Sigma_{i+1}, \ldots, \Sigma_{i+n}$, denoted by $T^{1/n}_{i+1, \ldots, i+n}$, as follows. Let $c$ be a curve surrounding 
the $n$ subsurfaces as shown in Figure \ref{Fig:FractionTwistSubSurfaces}. Then, the twist map is the homeomorphism 
with support on the surface separated by $c$ that sends $\Sigma_j$ to $\Sigma_{j+1}$
for  $j=i+1, \dots, (n+i-1)$ and sends $\Sigma_{i+n}$ to $\Sigma_{i+1}$ such 
that $(T^{1/n}_{{i+1}, \ldots, {i+n}})^n$ is a Dehn twist around the curve $c$. 
We then define
\[
\pi_a = \prod_{i=1}^{z(a)} 
T_{{z_a(i)}, \ldots, {\abs{a}+i}}^{1/(\abs{a} + i - z_a(i))}.
\]

Observe that, the $i$--th twist in the the above product moves the 
$\Sigma_{|a|+i}$ to $\Sigma_{z_a(i)}$. Therefore, the surfaces 
$\Sigma_{\abs{a}+1}, \dots, \Sigma_{\abs{a}+ z(a)}$ are sent by $\pi_a$ 
to surfaces in correspondence with the positions where $a$ has zeros. 
We finally define:
\begin{equation}\label{PhiMapAB}
    \begin{aligned}
    \Phi \from \Qinf &\to \FMap(\Sigma)\\
    a &\mapsto \pi_{a} h_{\sigma}^{\abs{a}}.
    \end{aligned}
\end{equation}

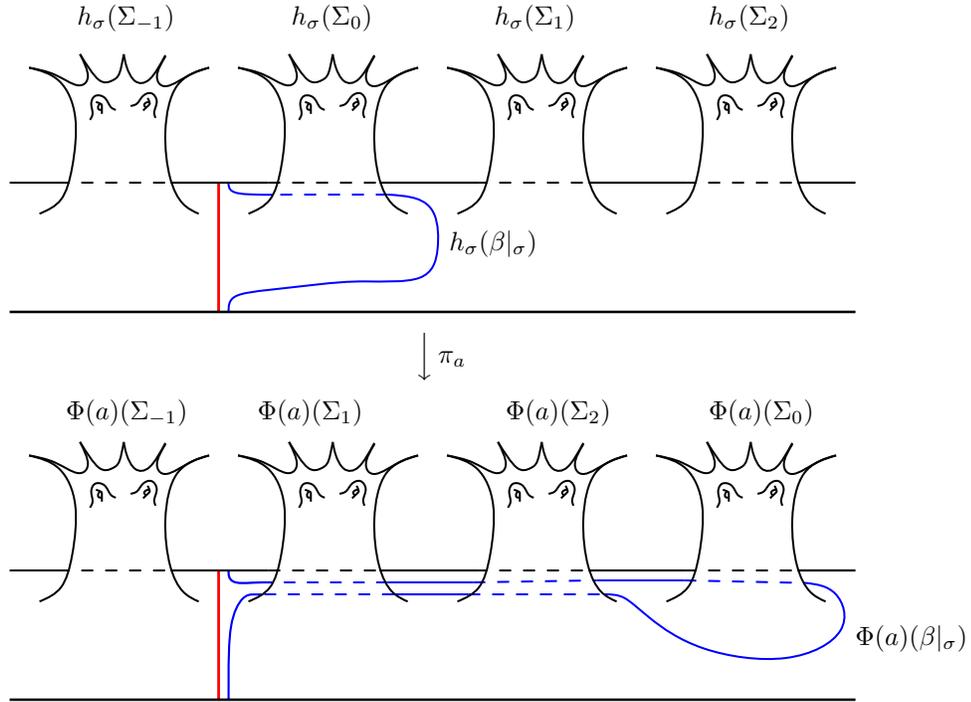
\begin{figure}[ht]
\setlength{\unitlength}{0.01\linewidth}
\begin{picture}(75,63)
\put(0,0){\includegraphics{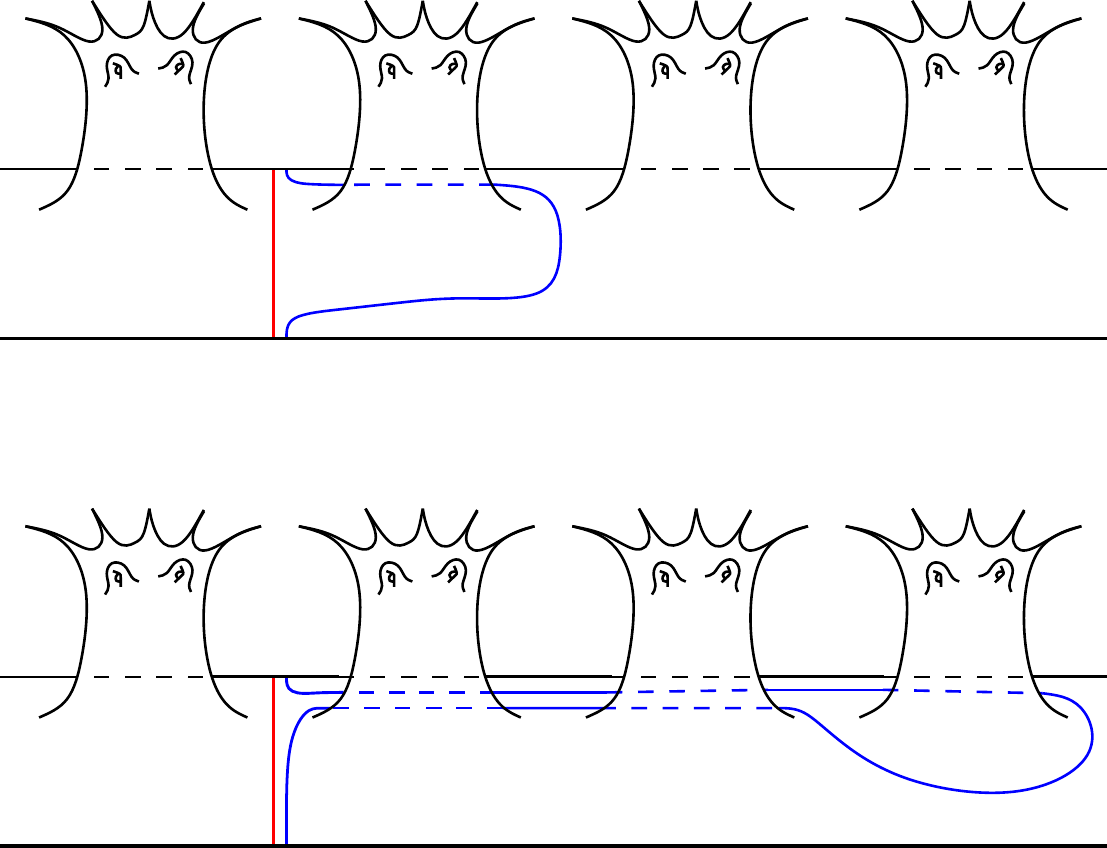}}
\put(36,30){$\Big\downarrow$}
\put(38,30){$\pi_a$}
\put(6,60){$h_\sigma(\Sigma_{-1})$}
\put(25,60){$h_\sigma(\Sigma_0)$}
\put(43,60){$h_\sigma(\Sigma_1)$}
\put(62,60){$h_\sigma(\Sigma_2)$}

\put(5,25){$\Phi(a)(\Sigma_{-1})$}
\put(22,25){$\Phi(a)(\Sigma_1)$}
\put(44,25){$\Phi(a)(\Sigma_2)$}
\put(62,25){$\Phi(a)(\Sigma_0)$}

\put(39,40){$h_\sigma(\beta|_\sigma)$}
\put(75,5){$\Phi(a)(\beta|_\sigma)$}

\end{picture} 
\caption{Consider the element $a = (0,0,1,0,\ldots) \in \Qinf$. 
Then $|a|=1$, $z(a)=2$ and hence $\pi_a =T^{1/2}_{2,3} \circ T^{1/2}_{1,2}$. 
The map $\Phi(a)= \pi_a h_\sigma$ maps the surface $\Sigma_0$ 
to $\Sigma_3$, There is only one index ($i=3$) where $a_i=1$ and hence 
one surface from $\Sigma_-$ is mapped to $\Sigma_+$ and is place at 
subsurface $\Sigma_i$ associated to this index.}
\label{Fig:WABPermutationExample} 
\end{figure}

Using the same reasoning as in Section \ref{Sec:SharkTank}, we can see that
$\Phi(a)$ sends some subsurface $\Sigma_j$ from $\Sigma_-$ to 
$\Sigma_+$ for every index $i>0$ where $a_i=1$ and sends this subsurface to $\Sigma_i$. If $b_i=1$ then
$\Sigma_i$ is sent back to $\Sigma_-$ under $\Phi(b)^{-1}$, otherwise
it stays in $\Sigma_+$. If $a_i=0$ the suburface sent to $\Sigma_i$
by $\Phi(a)$ comes from $\Sigma_+$. Now if $b_i=1$, then $\Sigma_i$
is sent to $\Sigma_-$ under $\Phi(b)^{-1}$. That is, the number
of surfaces that are sent from $\Sigma_\pm$ to $\Sigma_\mp$
under $\Phi(b)^{-1}\Phi(a)$ is exactly the number of indices $i$
where $a$ and $b$ are different, which is $|a-b|$. 

When $E^G$ is two sided, the surfaces $\Sigma_i$
are genus one surfaces with one boundary component. Therefore, 
sending one of these surfaces from $\Sigma_\pm$ to $\Sigma_\mp$ contributes 
2 to the norm defined in \eqnref{Eq:homology-norm} and we have
\begin{equation}\label{Eqn:QinfNormRelPhiEssentialShift-genus}
\norm{\Phi(b)^{-1}\Phi(a)(\beta)} = 2|a-b|. 
\end{equation} 
In the case where $E(z)$ is two-sided, $\Sigma_i$ has one maximal end which
is a point in $E(z)$. Therefore, 
sending one of these surfaces from $\Sigma_\pm$ to $\Sigma_\mp$ contributes 
1 to the norm defined in \eqnref{Eq:end-norm} and we have
\begin{equation}\label{Eqn:QinfNormRelPhiEssentialShift-end}
\norm{\Phi(b)^{-1}\Phi(a)(\beta)} = |a-b|. 
\end{equation} 
We now show that $\Phi$ is a quasi-isometric embedding from $\Qinf$ and 
$\FMap (\Sigma)$. 

\begin{proposition}\label{Prop:QIsometryQinfPhiAB}
Assume $E = E(\Sigma)$ is two-sided and let $\Phi \from \Qinf \to \FMap ( \Sigma)$
be the map from Equation \ref{PhiMapAB}. Then there exists $C>0$ 
such that, for every $a,b \in \Qinf$ we have:
\[
\frac 1C \cdot \abs{a-b} 
\leq \norm{\Phi(a) \Phi(b)^{-1}}_{\calS} \leq  C \cdot \abs{a-b}+3
\]
and hence $\Phi$ is a quasi-isometric embedding from 
$(Q^{\infty}, \abs{\param})$ to $(\FMap(\Sigma), \norm{\param}_{\calS})$.
\end{proposition}

\begin{proof} To obtain inequality on the left, we assume $C$ is larger
than the constants in Theorems \ref{Thm:LRelationToWordMetricEnds}
and \ref{Thm:HandleShiftBoundedNorm}. We then combine
\eqnref{Eqn:QinfNormRelPhiEssentialShift-end} with 
\thmref{Thm:HandleShiftBoundedNorm} 
in the case where $E^G$ is two-sided
and we combine \eqnref{Eqn:QinfNormRelPhiEssentialShift-genus}
with \thmref{Thm:LRelationToWordMetricEnds} 
in the case where $E(z)$ is 
two-sided. 

Recall that $x,y \in E$ are the ends of $\Sigma$ associated to $\sigma$. 
Let $A \in \calA$ be the set containing $x$ and $B \in \calA$ be the
set containing $y$. We now argue as in the proof of Theorem \ref{Prop:QIsometryQinfPhiSharkTank}, but in place 
punctures, we use disjoint subsurfaces $\Sigma_i$. We find homeomorphisms
$s_1, s_3 \in \nu_B$ and $s_2 \in \nu_A$ such that 
\[
(s_3 h_\sigma^{k_2} s_2 h_\sigma^{k_1} s_1) \circ 
\Phi(a) \Phi(b)^{-1} = id, 
\]
where $k_1$ is the number indices $i$ where $a_i=1$ and $b_i =0$ and
$k_2$ is the number of indices where $a_i=0$ and $b_i=1$.
Further assuming that $C \geq \norm{h_\sigma}_\calS$, we have 
\[
\norm{\Phi(a) \Phi(b)^{-1}}_\calS \leq C \cdot(k_1 + k_2) + 3 
= C \cdot |a-b| + 3 
\]
since the $s_i$ are generators.
This finishes the proof. 
\end{proof}

\thmref{Thm:Main}, restated below, now follows immediately.  

\begin{theorem}\label{Thm:FMapInfiniteAsympDim}
Assume $\Sigma$ has a two-sided end space. Then 
$\FMap(\Sigma)$ and $\Map(\Sigma)$ have infinite asymptotic dimension.
\end{theorem}

\begin{proof}
By Proposition \ref{Prop:QIsometryQinfPhiAB}, the map $\Phi \from (\Qinf, \abs{\param}) \to (\FMap ( \Sigma), d_\calS)$ 
is a quasi-isometric embedding and $\Qinf$ has infinite asymptotic dimension by 
\thmref{Thm:Infinite-Asymptotic-Dimension-Qinf}. Therefore $\FMap(\Sigma)$ also has infinite asymptotic dimension.

By \thmref{Thm:CompatibleGenSetQI}, $\FMap(\Sigma)$ and 
$\Map(\Sigma)$ are quasi-isometric. Fact $2.2$ tells us that 
asymptotic dimension is preserved under quasi-isometries, and hence 
$\Map(\Sigma)$ has infinite asymptotic dimension.
\end{proof}

\section{Equivalence of Algebraically and Topologically Essential}\label{Sec:EquivalenceAlgebraicTopologicallyEssential}

In this section we prove Theorems \ref{Thm:Existence} and \ref{Thm:Characterization}
from the introduction. Let $h_\sigma$ be a shift map as described in the introduction. Recall that the support
of $h_\sigma$ is a strip $\sigma$, containing a collection of sub-surfaces $\Sigma_i$, exiting towards points 
$x,y \in E$. 
%
We make use of the following theorem of Rosendal giving an equivalent condition for the notion of 
coarsely bounded sets that is more suitable for our purposes. 

\begin{theorem}[Rosendal, Prop. 2.7 (5) in \cite{Rose}] \label{Thm:Rosendall}
Let $A$ be a subset of a Polish group $G$. The following are equivalent
\begin{itemize}
    \item $A$ is coarsely bounded.
   \item For every neighborhood $\nu$ of the identity in $G$, there is a finite subset $\calF$ 
   and some $k \geq  1$ such that $A \subset (\calF\nu)^k$. 
\end{itemize}
\end{theorem} 

We prove several special cases of \thmref{Thm:Existence}. We then combine them to give a proof of the general case.  

\begin{proposition} \label{Prop:finite-genus}
Assume $\Sigma_i$ is a surface of genus $g$ with one boundary component. 
If $h_\sigma$ is essential, then $E^G$ is two-sided giving a decomposition $E^G = X \sqcup Y$
such that $x \in X$ and $y \in Y$. In particular, if $E^G$ is not two-sided, then $h_\sigma$
is not essential. 
\end{proposition} 

\begin{proof}
Let $A, B \in \calA$ be such that $x \in A$ and $y \in B$. There is a homeomorphic copy $\sigma'$
of $\sigma$ which is the concatenation of an infinite genus half-strip in $\Sigma_A$, 
a zero genus compact strip in $L$, and an infinite genus half-strip in $\Sigma_B$. 
Let $g$ be a homeomorphism sending $\sigma$ to $\sigma'$. For $\calF$ and $\nu$ as in 
\thmref{Thm:Rosendall}, let $\calF' = g^{-1} \calF g$ and $\nu' = g^{-1} \nu g$. Then 
\[
h_\sigma^n \in (\calF \nu)^k
\qquad \Longleftrightarrow  \qquad
h_{\sigma'}^n \in (\calF' \nu')^k.
\]
That is, $\langle h_\sigma \rangle$ is CB if and only if $\langle h_{\sigma'} \rangle$ is CB. 
In fact, in general, any set conjugate to a CB set is CB. Hence, it is enough to prove the Proposition for $h_{\sigma'}$.

\subsection*{Claim:} Assume there is no decomposition $E^G = X \sqcup Y$ where 
$x\in X$, $y \in Y$ and $X,Y$ are closed $\FMap(\Sigma)$--invariant sets. 
Then, there is a sequence 
\[
A = A_0, A_1, \dots , A_k = B, \qquad A_i \in \calA
\]
and a sequence of ends $z_i \in E^G$ and $z_i' \in E(z_i)$, $i = 0, \dots, (k-1)$ such that 
$z_0 \in A_0-\calM(A_0)$, $z_{k-1}' \in A_k-\calM(A_k)$ and, for $i = 1, \dots, (k-1)$ , 
$z_{i-1}', z_i \in A_i-\calM(A_{i-1})$. 

To see that this claim holds, let $\calA'$ be the subset of $\calA$ such that, for $C \in \calA'$, there is a sequence 
as described above starting from $A$ and ending in $C$. Let $\calB' = \calA - \calA'$. 
Then, for every $A' \in \calA'$, $B' \in \calB'$ and $z \in E_{A'}\cap E^G$ we have 
$E(z) \cap E_{B'} = \emptyset$. Define 
\[
X = \bigcup_{A' \in \calA'} E_{A'} \cap E^G
\qquad\text{and}\qquad
Y = \bigcup_{B' \in \calB' } E_{B'}\cap E^G.
\]
The sets $X$ and $Y$ are closed and $\FMap(\Sigma)$ invariant, and the definition of $X$ implies that 
we have $x \in X$. The assumption implies that $y$ cannot be in $Y$. Therefore $B \in \calA'$ and the 
above sequence exists, which proves the claim.

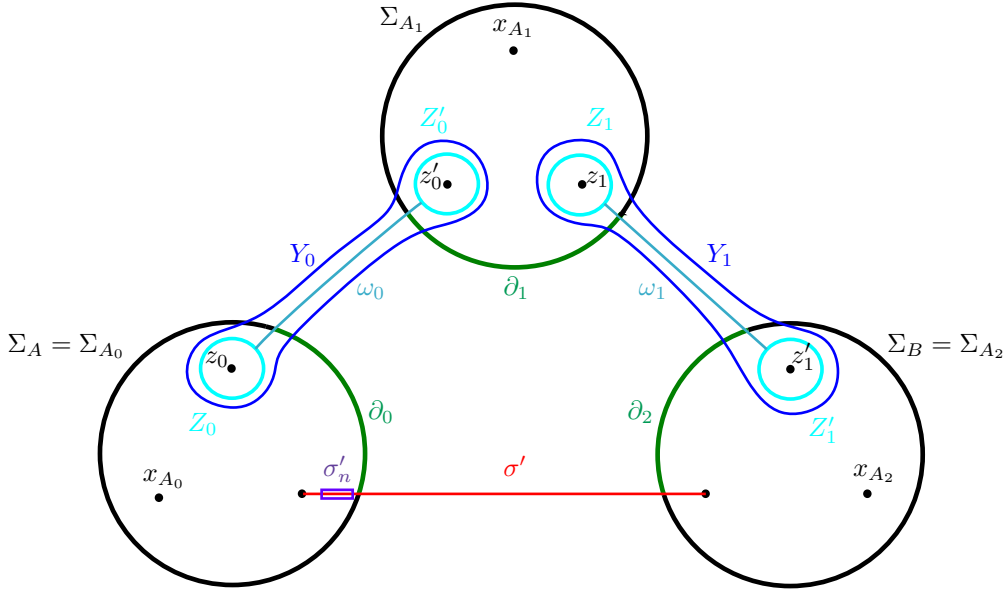
\begin{figure}[ht]
\setlength{\unitlength}{0.01\linewidth}
\begin{picture}(75,53)
\put(0,0){\includegraphics[width=110mm]{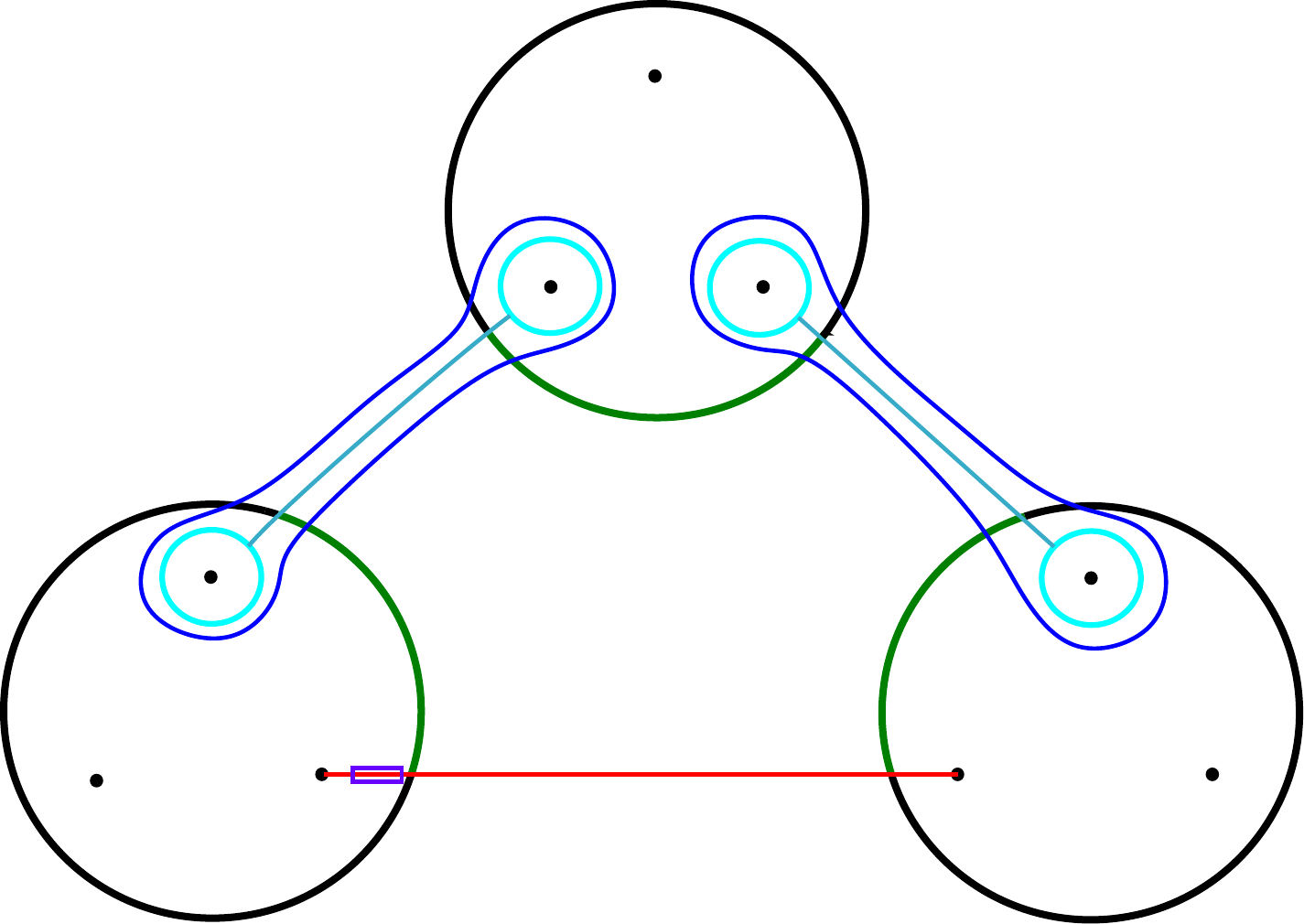}}
\put(-8,21){$\Sigma_A=\Sigma_{A_0}$}
\put(25,50){$\Sigma_{A_1}$}
\put(70,21){$\Sigma_B=\Sigma_{A_2}$}
\put(4,9.5){$x_{A_0}$}
\put(35,49.5){$x_{A_1}$}
\put(67,10){$x_{A_2}$}
\put(9.5,20){$z_{0}$}
\put(8,14){\textcolor{cyan1}{$Z_{0}$}}
\put(17,29){\textcolor{blue}{$Y_{0}$}}
\put(23,26){\textcolor{dustyblue}{$\omega_{0}$}}
\put(54,29){\textcolor{blue}{$Y_{1}$}}
\put(48,26){\textcolor{dustyblue}{$\omega_{1}$}}
\put(28.5,36){$z_{0}'$}
\put(43.25,36){$z_{1}$}
\put(28.5,41){\textcolor{cyan1}{$Z_{0}'$}}
\put(43.25,41){\textcolor{cyan1}{$Z_{1}$}}
\put(61.5,20){$z_{1}'$}
\put(63,13.5){\textcolor{cyan1}{$Z_{1}'$}}
\put(24,15){\textcolor{ForestGreen}{$\partial_0$}}
\put(36,26){\textcolor{ForestGreen}{$\partial_1$}}
\put(47,15){\textcolor{ForestGreen}{$\partial_2$}}
\put(36,10){\textcolor{red}{$\sigma'$}}
\put(20,10){\textcolor{RoyalPurple}{$\sigma_n '$}}
\end{picture}
\caption{This is the picture of the our set up for $k=2$. The arcs $\omega_i$ and $\partial_i$ are chosen 
so that, the concatenation of $\partial_0, \omega_0|_L, \partial_1, \omega_1|_L$ is homotopic to 
the concatenation of $\sigma'|_L$ followed by $\partial_2$.} 
\label{Fig:not-essential} 
\end{figure}

For $i=0, \dots, (k-1)$, fix disjoint sub-surfaces $Z_i \subset \Sigma_{A_i} -x_{A_i}$ and 
$Z_i'\subset \Sigma_{A_{i+1}}-x_{A_{i+1}}$ such that $z_i$ is an end of $Z_i$, $z_i'$ is an end of $Z_i'$ and 
$Z_i$ is homeomorphic to $Z_i'$, for example, we can choose $Z_i$ so that $E(Z_i)$ is an stable neighborhood
of $z_i$. Connect $Z_i$ to $Z_i'$ by an arc $\omega_i$ that is contained in 
$L \cup \Sigma_{A_i} \cup \Sigma_{A_{i+1}}$. Let $Y_i$ be a regular neighborhood of 
$Z_i \cup Z_i' \cup \omega_i$. Further, we assume that subsurfaces $Y_i$ are disjoint from each other 
and are disjoint from the strip $\sigma'$ (see \figref{Fig:not-essential}). For $i =1, \dots (k-1)$, let $\partial_i$
be a sub-arc of $\partial \Sigma_{A_i}$ connecting $\omega_{i-1} \cap \partial \Sigma_{A_i}$ to 
$\omega_{i} \cap \partial \Sigma_{A_i}$. Also, let $\omega_i|_L$ be the restriction 
of $\omega_i$ to $L$. We further assume that the concatenation of $\omega_i|_L$ and $\partial_i$, 
which is an arc connecting $\partial \Sigma_A$ to $\partial \Sigma_B$, is homotopic (relative
$\partial \Sigma_A \cup \partial \Sigma_B$) to either arc in the the restriction of $\partial \sigma'$ to 
$L$, which we denote by $\partial \sigma'|_L$. Let $\partial_0$ be the sub-arc of $\partial \Sigma_A$
connecting a point in $\partial \sigma' \cap \partial \Sigma_A$ to the points $\omega_0 \cap \partial \Sigma_A$
and similarly let $\partial_k$ be the sub-arc of $\partial \Sigma_B$ connecting a point in 
$\partial \sigma' \cap \partial \Sigma_B$ to the points $\omega_{k-1} \cap \partial \Sigma_B$. 
Again, these can be chosen such that the concatenation $(\partial_i \omega_i|_L )_{i=0}^{(k-1)}$ 
is homotopic to the concatenation of an arc in $\partial \sigma'|_L$ and $\partial_k$ relative their end points. 

Let $g_i$ be a homeomorphism that sends $Z_i$ to $Z_i'$, whose support 
is in $Y_i$ and that sends $\omega_i$ to itself in the reverse direction. 
Let $\calF = \{g_i\}_{i=0}^{(k-1)}$. 

We now show that, for every $n \in \mathbb Z$, $h_{\sigma'}^n \in (\nu_L \calF)^{k+1}$. Let $\sigma'_n$ be 
a sub-strip of $\sigma'$ of genus $n$ so that the genus of what remains between $\sigma_n'$ and
$\partial A$ is zero. Choose an element $f_0 \in \nu_{A_0}$ that sends $\sigma'_n$ to a 
$Z_0$, which is possible since $z_0 \in E^G$ and therefore  $Z_0$ has infinite genus. In fact, we can do this in a way such that 
$f_0(\sigma')\cap \Sigma_A$ starts the same as $\sigma'\cap\Sigma_A$ away from a small neighborhood of 
$\partial \Sigma_A$, then it follows $\partial_0$ (staying in a small neighborhood of $\partial \Sigma_A$) 
to the intersection point of $\partial \Sigma_A$ and $\omega_0 \cap \Sigma_A$, then follows 
$\omega_0 \cap \Sigma_A$ (staying in $Y_0$) to $Z_0$, then comes back the same way along 
$\omega_0 \cap\Sigma_A$ and $\partial_0$ and continues along $\sigma'$ into $B$
(see the left picture in \figref{Fig:Steps}). Then, $g_0 \circ f_0(\sigma_n') \subset Z_0'$. 

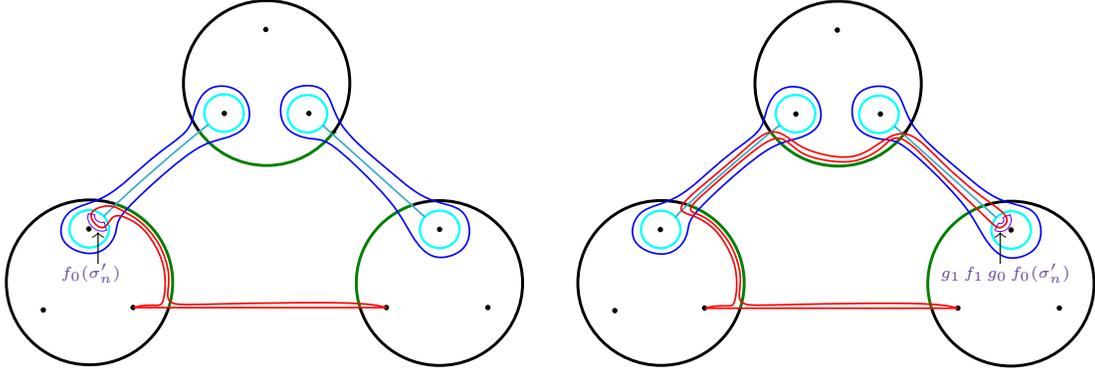
\begin{figure}[ht]
\setlength{\unitlength}{0.01\linewidth}
\begin{picture}(75,30)
\put(-10,0){\includegraphics[width=145mm]{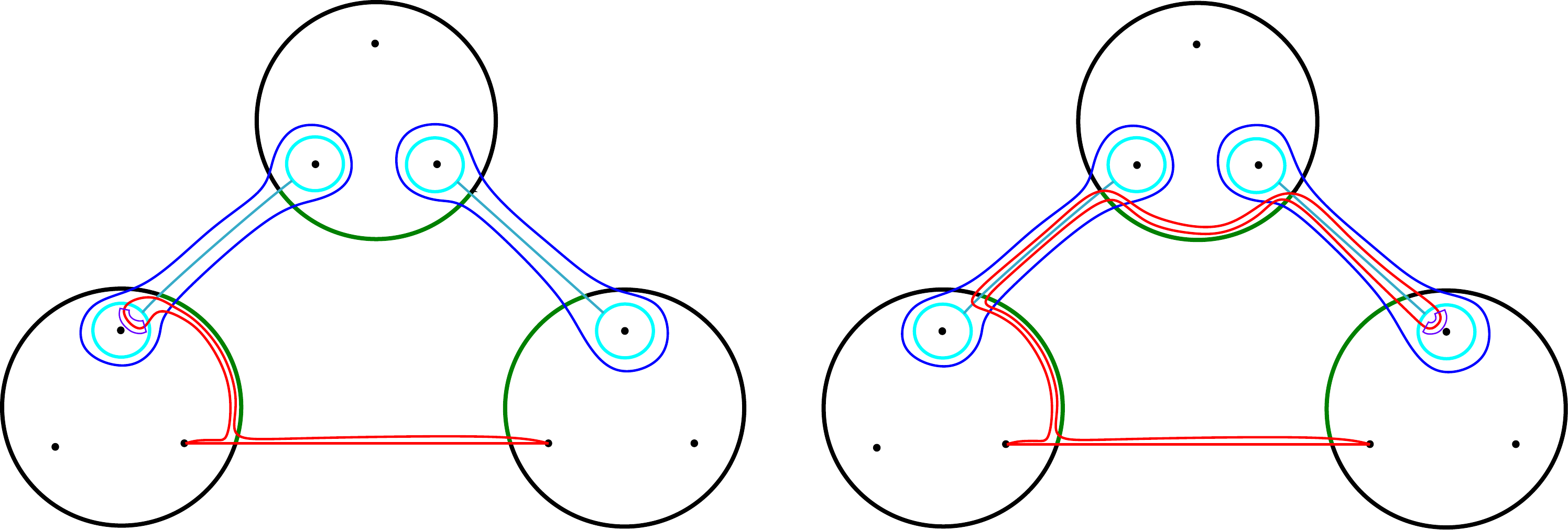}}
\put(-5,8){\textcolor{RoyalPurple}{\tiny{$f_0(\sigma_n')$}}}
\put(-2.5,10){$\big\uparrow$}
\put(73,8){\textcolor{RoyalPurple}{\tiny{$g_1 \, f_1 \, g_0 \, f_0(\sigma_n')$}}}
\put(77.5,10){$\big\uparrow$}
\end{picture}
\caption{The picture on the left depicts the image of $\sigma'$ and $\sigma'_n$ under the map 
$f_0$. The picture on the right depicts the image of $\sigma'$ and $\sigma'_n$ under the map 
$\prod_{i=0}^{k-1} g_i \circ f_i (\sigma')$. The strip $\prod_{i=0}^{k-1} g_i \circ f_i (\sigma')$ can be 
homotoped to $\sigma'$ inducing the shift map $h_{\sigma'}^n$. }
\label{Fig:Steps} 
\end{figure}

We then find a $f_1 \in \nu_{A_1}$ that sends $g_0 \circ f_0(\sigma_n')$ to $Z_1$. We can do this in a way such that
$f_1 \circ g_0 \circ f_0 (\sigma')$ starts the same as $\sigma'$ inside $\Sigma_A$ (away from a small neighborhood 
of $\partial \Sigma_A$) then it follows $\partial_0$, then $\omega_0|_L$, then $\partial_1$ to the intersection 
point of $\partial \Sigma_{A_1}$ and $\omega_1$, then it follows $\omega_1 \cap \Sigma_{A_1}$ 
(staying in $Y_1$) to $Z_1$, then comes back the same way along $\omega_1 \cap \Sigma_{A_1}$, $\partial_1$, 
$\omega_0|_L$ and $\partial_0$, and then follows $\sigma'$ into $B$. We then apply $g_1$ and we have 
$g_1 \circ f_1 \circ g_0 \circ f_0(\sigma_n') \subset Z_1'$. 

Following the same argument, we can find $f_i$, $i =0, \dots, (k-1)$ such that 
\[
\left( \prod_{i=0}^{k-1} g_i \circ f_i\right) (\sigma_n') \subset Z_{k-1}' \subset \Sigma_{A_k} = \Sigma_B.
\]
The strip $\prod g_i \circ f_i(\sigma')$ starts the same as $\sigma'$ in $\Sigma_A$, then 
follows the concatenation of the arcs $\omega_i|_L$ and $\partial_i$ (which by assumption is homotopic to arcs 
$\partial \sigma'|_L$), then follows $\omega_k\cap \Sigma_B$ to $Z_k$, then back the same way to 
$\partial \Sigma_A \cap \sigma'$, and then continues the same as $\sigma'$ (see the right picture in 
\figref{Fig:Steps}). 

The portion of $\prod_{i=0}^{k-1} g_i \circ f_i (\sigma_n')$ traveling back from $\partial \Sigma_B$
to $\partial \Sigma_A$ has genus zero and can be homotoped into $\Sigma_B$. The portion going forward from 
$\partial \Sigma_A$ to $\partial \Sigma_B$ also has genus zero and can be homotoped to a neighborhood of 
the concatenation of $\partial \sigma'|_L$ and $\partial_k$. Hence, there is $f_k \in \nu_B$ such that 
$f_k \circ \prod g_i \circ f_i (\sigma')$ can be homotoped into $\sigma'$ and hence 
$f_k \circ \prod g_i \circ f_i $ can be considered as a map with support in $\sigma'$. 
Since, $\sigma_n'$ has moved from $A$ to $B$, we in fact have $f_k \circ \prod g_i \circ f_i= h_{\sigma'}^n$. 
That is $h_{\sigma'}^n \in (\nu_L \calF)^{k+1}$. This finishes the proof. 
\end{proof}

\begin{proposition} \label{Prop:single-end}
Assume $\Sigma_i$ is a surface with either zero or infinite genus and where $E(\Sigma_i)$ has a single maximal 
point $z$. If $h_\sigma$ is essential, then $E(z)$ is two-sided giving the decomposition 
$\Accu(x) = X \sqcup Y$ where $x \in X$ and $y \in Y$. 
\end{proposition} 

\begin{proof}
The proof is nearly identical to the proof of \propref{Prop:finite-genus}. We outline it here. 
Let $x$, $y$ be the ends of $\Sigma$ associated to the strip $\sigma$. Let $A, B \in \calA$ be
such that $x \in A$ and $y \in B$. First find $g \in \FMap(\Sigma)$ such that $\sigma'=g(\sigma)$ 
starts in $\Sigma_A$, continues in $L$, and then enters $\Sigma_B$. It is enough to to show that 
the group generated by $h_{\sigma'} = g h_\sigma g^{-1}$ is CB. 

Note that $x,y \in \Accu(z)$. We can find a sequence
\[
A = A_0, A_1, \dots , A_k = B, \qquad A_i \in \calA
\]
and a sequence of ends $z_i \in \Accu(z)$ and $z_i' \in E(z_i)$, $i = 0, \dots, (k-1)$ such that 
$z_0 \in A_0-\calM(A_0)$, $z_{k-1}' \in A_k-\calM(A_k)$ and, for $i = 1, \dots, (k-1)$ , 
$z_{i-1}', z_i \in A_i-\calM(A_{i-1})$. A similar argument to the proof of \propref{Prop:finite-genus}
shows that if such a sequence does not exist, then $E(z)$ is two-sided which would be a contradiction. 

We then construct sub-surfaces $Z_i$ and $Y_i$ and maps $g_i$ as before, and set 
$\calF = \{ g_i\}_{i=0}^{k-1}$. Let $\sigma_n'$ be a sub-strip of $\sigma'$ that contains 
the $n$ sub-surfaces $\Sigma_i$ that are nearest to $\partial \Sigma_A$. As in the 
proof of \propref{Prop:finite-genus}, we can choose  $f_i \in \nu_{A_i}$ such that $f_k \circ \prod (g_i \circ f_i)$ 
sends $\sigma'$ to the strip that is homotopic to $\sigma'$ and moves $\sigma_n$ from $\Sigma_A$ to
$\Sigma_B$. That is 
\[
h_{\sigma'}^n = f_k \circ \prod (g_i \circ f_i) \in \big(\nu_L \calF \big)^{k+1}.
\] 
This finishes the proof.
\end{proof}

\begin{proposition} \label{Prop:cantor-set}
Assume $\Sigma_i$ is a surface with either zero or infinite genus, and where $\calM(E(\Sigma_i))$ is a Cantor set. 
Then $h_\sigma$ is not essential. 
\end{proposition}

\begin{proof}
As before, after conjugation, we can assume $\sigma$ intersects only $\Sigma_A \cup L \cup \Sigma_B$,  where
$\Sigma_i \in \Sigma_A$ for $i \leq 0$ and $\Sigma_i \in \Sigma_B$ for $i>0$. 
Consider a simple closed curve in the strip $\sigma$ separating $\Sigma_i$ and $\Sigma_{i+1}$ 
from the rest of the strip. Denote this subsurface by $\Sigma_{i, i+1}$. If the set of maximal 
ends in $\Sigma_i$ is a Cantor set, the end space of $\Sigma_i$ is homemorphic to the end
space of $\Sigma_{i, i+1}$ (see \figref{Fig:Cantor-End}).

\begin{figure}[h!]
    \centering
    \includegraphics{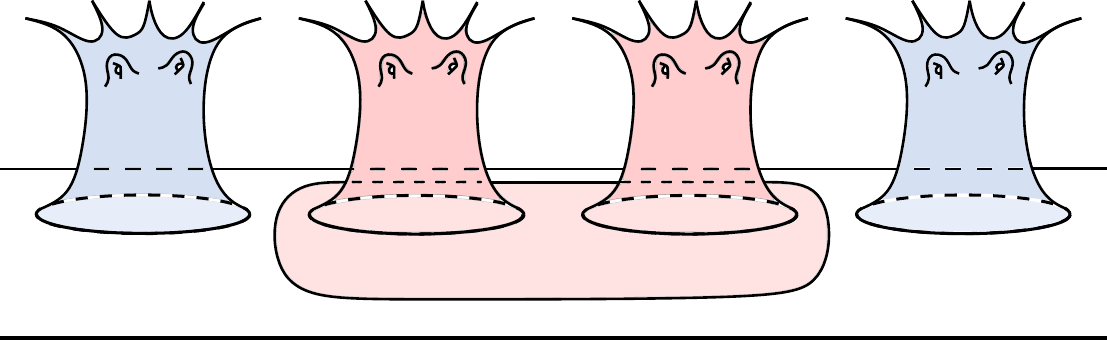}
    \put(-169,20){$\Sigma_{i,i+1}$}
    \put(-287,60){$\Sigma_{i-1}$}
    \put(-51,60){$\Sigma_{i+2}$}
    \label{Fig:Cantor-End}
    \caption{When $\calM(E(\Sigma_i))$ is a cantor set the surface $\Sigma_{i, i+1}$ is homeomorphic to 
    $\Sigma_i$.}
\end{figure}

Similarly, the subsurface $\Sigma_{i, i+n}$, which is a 
subsurface
of $\sigma$ with one boundary component containing the sub-surfaces $\Sigma_i, \Sigma_{i+1}, \dots,
\Sigma_{i+n}$, is also homeomorphic to $\Sigma_i$. Hence, for every $n$, there is a homeomorphism 
$f_n \in \nu_A$ with support in $\sigma$ such that $\Sigma_i$  is sent to $\Sigma_{i+n}$ for 
$i < -n$, $\Sigma_i$ is sent to itself for $i>0$, and the surface $\Sigma_{-n, 0}$ is sent to 
$\Sigma_0$. Similarly, there  there is a homeomorphism 
$g_n \in \nu_A$ with support in $\sigma$ such that $\Sigma_i$  is sent to $\Sigma_{i+n}$ for 
$i >1$, $\Sigma_i$ is sent to itself for $i\leq 0$, and the surface $\Sigma_{1}$ is sent to 
$\Sigma_{1,n+1}$. Then $h_\sigma^n = g_n h_\sigma f_n$. That is, we can collect the surfaces 
$\Sigma_{-n} \dots \Sigma_0$ to the surface $\Sigma_0$, shift $\sigma$ once and the expand
these surfaces to $\Sigma_1, \dots, \Sigma_{n+1}$, which means we have written $h_\sigma^n$ as a composition of $3$ 
homeomorphisms. By letting $\calF = \{ h_\sigma\}$, we have that
$\langle h_\sigma \rangle$ is contained in $(\nu_L \calF)^2$.  Therefore, by \thmref{Thm:Rosendall}, 
$\langle h_\sigma \rangle$ is a CB subset of $\FMap(\Sigma)$ and $h_\sigma$ is not essential. 
\end{proof}

We are now in a position where we can prove Theorem \ref{Thm:Existence}, which states that $\Map(\Sigma)$ 
contains an essential shift if and only if the endspace of $\Sigma$ is two-sided.

\begin{proof}[Proof of \thmref{Thm:Existence}]
One direction is already proven by \propref{Prop:2-sided-implies-essential}. 
It remains to show that, if $E$ is not two-sided then there is no essential shift map. 

Assume $E^G$ is not two-sided and there does not exist any $z \in E$ such that $E(z)$ is two-sided.
Let $h_\sigma$ be any shift map with support on a strip $\sigma$, containing the subsurfaces $\Sigma_i$, 
exiting towards $x, y \in E$. The end space $E(\Sigma_i)$ has finitely many maximal types each of which is either a 
Cantor set or a finite set. Therefore we can find finitely many sub-surfaces $\Sigma_i^j \subset \Sigma_i$, 
$j=0 , \dots, \ell$, each with one boundary component such that:
\begin{enumerate}
\item The surfaces $\Sigma_i^j$ are disjoint. Furthermore, $\Sigma - \Sigma_i^j$ is a compact planar surface
which means $E(\Sigma_i) = \sqcup_j E(\Sigma_i^j)$. 
\item For $j>0$, $\Sigma_i^j$ has zero genus or infinite genus. That is, if $\Sigma_i$ has finite genus
we include all the genus in the surface $\Sigma_i^0$ so every other subsurface of $\Sigma_i$ contains no genus. If
$\Sigma_i$ already has zero genus or infinite genus, then there is no need to choose $\Sigma_i^0$. We can choose
$\Sigma_i^0$ to be a disk or have $j$ range from $1$ to $\ell$. Also, if $E(\Sigma_i^j) \cap E^G =\emptyset$
then we make sure $\Sigma_i^j$ has genus zero. 
\item For every $j=0, \dots, \ell$, the surfaces $\Sigma_i^j$ are homeomorphic for all $i \in \ZZ$. 
That is we decompose each $\Sigma_i$ in the same way. 
\end{enumerate}
Now the strip $\sigma$ can be decomposed to parallel strips $\sigma^j$ each containing sub-surfaces 
$\Sigma_i^j$. Since $h_{\sigma^j}$ have disjoint support, they commute and $h_\sigma = \prod h_{\sigma^j}$.
That is, the group generated by $h_{\sigma_j}$, $j=0, \dots, \ell$, is an abelian group that contains 
the group generated by $h_\sigma$. Hence, if each $\langle h_{\sigma^j} \rangle$ is CB then 
$\langle h_{\sigma^0, \dots, h_{\sigma^\ell}} \rangle$ is also CB and thus $\langle h_{\sigma} \rangle$ 
is CB. 

Since $E^G$ is not two-sided, by \propref{Prop:finite-genus}, $h_{\sigma^0}$ is not essential and 
$\langle h_{\sigma^0} \rangle$ is CB. For $j=1, \dots, \ell$, if $\calM(E(\Sigma_i^j))$ is a single 
point then $\langle h_{\sigma^0} \rangle$ is CB by \propref{Prop:single-end}. If $\calM(E(\Sigma_i^j))$
is a Cantor set, then  $\langle h_{\sigma^0} \rangle$ is CB by \propref{Prop:cantor-set}. 
This finishes the proof. 
\end{proof}

Finally, we prove Theorem \ref{Thm:Characterization} which tells us that a shift map is essential 
if and only if there is either a decomposition of the ends accumulated by genus, or there is a 
decomposition of the accumulation set of a maximal point of the surface.

\begin{proof}[Proof of \thmref{Thm:Characterization}]
We begin by proving the forward direction.
Consider the decomposition of the shift map $h_\sigma = \prod_j h_{\sigma^j}$ constructed in the 
proof of Theorem \ref{Thm:Existence}. As mentioned before, 
if $h_\sigma$ is essential then some  $h_{\sigma_j}$ is essential. If $h_{\sigma_0}$ is essential, then since $\Sigma_i^0$ contains finite genus, 
\propref{Prop:finite-genus} implies the first bullet point in \thmref{Thm:Characterization} holds. If $j>0$ and 
$\calM(E(\Sigma_i^j))$ is a single point, then by \propref{Prop:single-end} the second bullet point in
\thmref{Thm:Characterization} holds. Finally, $\calM(E(\Sigma_i^j))$ cannot be a Cantor set
by \propref{Prop:cantor-set}, which proves the forward direction. 

In the other direction, first suppose that $E^G$ is two-sided giving a decomposition $E^G = X \sqcup Y$ 
where $x \in X$ and $y \in Y$, then the proof of \propref{Prop:2-sided-implies-essential}
shows that $h_{\sigma^0}$ is essential. In fact, there is a length function $\norm{\param}$
such that, for $n>0$, $\norm{h_{\sigma^0}^n} = 2n$. Since, in this case, $\Sigma_i^j$ have genus zero 
for $j>0$, the shift maps $h_{\sigma^j}$ act trivially on homology. Hence $\norm{h_{\sigma}^n} = 2n$
and hence $h_\sigma$ is essential. 

Similarly, if there exists some $z \in \calM(\Sigma_i)$ such that $E(z)$ is two-sided giving a 
decomposition $\mathrm{Accu}(z)= X \sqcup Y$ where $x \in X$ and $y \in Y$, then the proof of 
\propref{Prop:2-sided-implies-essential}
shows that $h_{\sigma^j}$ is essential for some $j>0$. In fact, there is a length function $\norm{\param}$
such that, for $n>0$, $\norm{h_{\sigma^0}^n} = n$. It is possible that the $\Sigma_i^j$ are homeomorphic 
for several different $j$. We can assume, after re-indexing, that $\Sigma_i^1, \dots, \Sigma_i^{\ell'}$
all have same type of countable maximal end. Then $\norm{h_{\sigma}^n} = \ell' \cdot n$
and hence $h_\sigma$ is essential. 
\end{proof}

\section{Two sources of non-trivial topology}\label{Sec:NonTrivialGeometry}

In this section, we prove \thmref{Thm:Non-Trivial-Geometry} which states that if $\Sigma$ does not
have a non-displaceable subsurface and $\FMap(\Sigma)$ does not have an essential shift then 
$\Map(\Sigma)$ and $\FMap(\Sigma)$ are CB.
This implies that $\Map(\Sigma)$ and $\FMap(\Sigma)$ are quasi-isometric to a point and therefore do not have an interesting
geometry. To prove this theorem, we make use of the notion of an avenue surface first introduced in \cite{HQR}. 
Recall that we always assume that $\Sigma$ is stable. 

\begin{definition}
An \textit{avenue surface} is a connected orientable surface $\Sigma$ which does not contain any 
non-displaceable subsurface of finite type, and whose mapping class
group $\Map(\Sigma)$ is CB-generated but not CB.
\end{definition}

That is, the only possible examples for surfaces that have no non-displaceable 
subsurfaces, no essential shifts, and are not CB are the avenue surfaces. 
A topological description of avenue surfaces was given in \cite{HQR}. 

\begin{lemma}[Lemmas 4.5 and 4.6 in \cite{HQR}] \label{Lem:Avenue}
Let $\Sigma$ be an avenue surface. Then $\Sigma$ has either 0 or infinite genus, and
$\Sigma$ has exactly two ends of maximal type, that is, $\calM(E)=\{ x_1, x_2\}$. 
Furthermore, for every $x \in E-\{x_1, x_2\}$, the set $E(x)$ accumulates to both $x_1$
and $x_2$.
\end{lemma}

In order to state the classification of CB mapping class groups from \cite{MR}, we must first 
define the notion of self-similarity for a space of ends, and the notion of telescoping for 
an infinite-type surface.

\begin{definition} \label{Def:SS-Telescoping}
A space of ends $(E, E^G)$ is said to be {\em self-similar} if for any decomposition 
$E = E_1 \sqcup E_2 \sqcup \ldots \sqcup E_n$ of $E$ into pairwise disjoint clopen sets, 
there exists a clopen set $D$ in some $E_i$ such that $(D, D\cap E^G)$ is homeomorphic to $(E, E^G)$.  

A surface $\Sigma$ is \textit{telescoping} if there are ends $x_1, x_2 \in E$ and disjoint 
clopen neighborhoods $V_i$ of $x_i$ in $\Sigma$ such that for all clopen neighborhoods 
$W_i \subset V_i$ of $x_i$, there exist homeomorphisms $f_i, h_i \in \FMap(\Sigma)$, with 
\[
f_i(W_i) \supset (\Sigma - V_{3-i}) 
\qquad h_i(W_i) = Vi, 
\qquad\text{and}\qquad 
h_i(V_{3-i}) = V_{3-i}.
\]
\end{definition}

A classification of CB mapping class groups was given in \cite{MR}. 

\begin{theorem}[Theorem 1.7 in \cite{MR}] \label{Thm:CB:Classification}
The group $\Map(\Sigma)$ is CB if and only if $\Sigma$ has infinite or zero genus and $E$
is either self-similar or telescoping. 
\end{theorem}

We are now ready to prove \thmref{Thm:Non-Trivial-Geometry} restated below. 
\newline

\noindent\textbf{Theorem \ref{Thm:Non-Trivial-Geometry}.} (Two sources of non-trivial geometry).
\textit{
If $\Map(\Sigma)$ does not have an essential shift and $\Sigma$ does not contain 
a non-displaceable subsurface then $\Map(\Sigma)$ is quasi-isometric to a point.}

\begin{proof}
By way of contradiction, suppose that $\Sigma$ does not have a non-displaceable subsurface, $\FMap(\Sigma)$ 
does not contain an essential shift, and that $\FMap(\Sigma)$ is not CB.
Then by assumption, $\Sigma$ is an avenue surface. 
By \lemref{Lem:Avenue}, 
$\calM(E)=\{x_1, x_2\}$ and for any other $x \in E - \{x_1, x_2\}$, $E(z)$ accumulates to 
both $x_1$ and $x_2$. Consider an end $z \in \calM(E-\{x_1, x_2\})$, that is, an end that is
maximal in $E-\{x_1, x_2\}$. If $E(z)$ is countable, then $E$ is two-sided with 
$X=\{x_1\}$ and $Y=\{x_2\}$. This cannot happen since we are
assuming there are no essential shifts (see \thmref{Thm:Existence}).
Therefore, $E(z)$ is a Cantor set. 

Furthermore, either $\Sigma$ has genus zero or it is infinite genus
and $\{x_1, x_2 \} \subset E^G$ (otherwise, $\Map(\Sigma)$ is not even locally CB, see 
\cite[Theorem 1.4]{MR}). We now show that these assumptions imply 
that $\Sigma$ is telescoping and, by \thmref{Thm:CB:Classification}, $\Map(\Sigma)$ is CB, which will prove our claim.

We check the definition of telescoping. For $i=1,2$, let $V_i$ be disjoint stable neighborhoods 
of $x_i$ such that $E-(V_1 \cup V_2)$ contains every maximal type in $E-\{x_1, x_2\}$.
Let $W_i$ be the given smaller neighborhoods of $x_i$.  Since the $V_i$ are stable neighborhoods, 
we have that $V_i$ is homeomorphic to $W_i$. We claim that $E-(V_1 \cup V_2)$ is homeomorphic to 
$E-(V_1 \cup W_2)$. Since all maximal types in $E-\{x_1, x_2\}$ are present in $E-(V_1 \cup V_2)$, 
the sets $E-(V_1 \cup V_2)$  and $E-(V_1 \cup W_2)$ have the same maximal types. But all
these types are Cantor sets. Therefore, for every $y \in V_2-W_2$ there is a $z \in E-(V_1 \cup V_2)$
such that $E(y)$ accumulates to $z$. By \cite[Lemma 4.18]{MR}, there are small neighborhoods $U_y$
of $y$ and $U_z$ of $z$ such that $U_y \cup U_z$ is homeomorphic to $U_z$. The neighborhoods
$U_y$ give a covering of $V_2-W_2$, hence there is a finite sub-covering. By making the neighborhoods
smaller, we can assume that they are disjoint. Since each neighborhood $U_y$ can be absorbed into $U_z$, 
the set $V_2-W_2$ can be absorbed into $E-(V_1 \cup V_2)$ and hence, $E-(V_1 \cup V_2)$ is
homeomorphic to $E-(V_1 \cup W_2)$. 

Now consider surfaces $\Sigma_{V_i}$ and $\Sigma_{W_i}$
whose end points are $V_i$ and $W_i$. We can choose these surfaces so that $\Sigma_{W_i} \subset \Sigma_{V_i}$
and so that $\Sigma_{V_i}$ is disjoint from $\Sigma_{V_{3-i}}$. Also, we can assume  
$\Sigma_{V_i} - \Sigma_{W_i}$ and $\Sigma-(\Sigma_{W_1} \cup \Sigma_{W_2})$ all have 
genus zero or infinity. The end space of $\Sigma_{V_i}- \Sigma_{W_i}$
is homeomorphic to the end space of $\Sigma-(\Sigma_{W_1} \cup \Sigma_{W_2})$
and they have the same 
genus. Also $\Sigma_{V_1}$ and $\Sigma_{W_1}$ have homeomorphic end spaces and the same
genus. Hence, there is a map $h_1$ such that $h_1(\Sigma_{W_1})=\Sigma_{V_1}$ and 
$h_1(\Sigma_{V_2}) = h_1(\Sigma_{V_2})$. The map $h_2$ in \defref{Def:SS-Telescoping} can 
also be similarly constructed.

The construction of the maps $f_i$ follows similarly. As above, we can send 
$\Sigma_{W_1}$ to $\Sigma-\Sigma_{V_2}$ while fixing some neighborhood
$W_2'$ of $x_2$ as long as $W_2'$ is small enough so that $V_2-W_2'$ intersect $E(z)$ for every 
$z \in \calM(E-\{x_1, x_2\})$. This proves that $\Sigma$
is telescoping and hence $\Map(\Sigma)$ is CB. Note that since
$\FMap(\Sigma)$ and $\Map(\Sigma)$ are quasi-isometric this means
$\FMap(\Sigma)$ is also CB which proves the theorem.
\end{proof}


\medskip

\bibliographystyle{plain}

\end{document}